\documentclass[11pt]{article}
\usepackage{amsfonts,amsmath,amssymb}
\usepackage{hyperref}
\usepackage{graphicx}
\usepackage{color}

\newcommand{\halmos}{\rule{1ex}{1.4ex}}
\def \qed {\nopagebreak{\hspace*{\fill}$\halmos$\medskip}}

\newtheorem{theorem}{Theorem}[section]
\newtheorem{proposition}[theorem]{Proposition}
\newtheorem{corollary}[theorem]{Corollary}
\newtheorem{conjecture}[theorem]{Conjecture}
\newtheorem{lemma}[theorem]{Lemma}

\newtheorem{remark}[theorem]{Remark}
\newtheorem{defi}[theorem]{Definition}

\newcommand{\bt}{\begin{theorem}}
\newcommand{\et}{\end{theorem}}
\newcommand{\bl}{\begin{lemma}}
\newcommand{\el}{\end{lemma}}
\newcommand{\bprop}{\begin{proposition}}
\newcommand{\eprop}{\end{proposition}}
\newcommand{\bcor}{\begin{corollary}}
\newcommand{\ecor}{\end{corollary}}
\newcommand{\br}{\begin{remark}\rm}
\newcommand{\er}{\end{remark}}
\newcommand{\bcon}{\begin{conjecture}}
\newcommand{\econ}{\end{conjecture}}
\newcommand{\bd}{\begin{defi}}
\newcommand{\ed}{\end{defi}}

\newcommand{\be}{\begin{equation}}
\newcommand{\ee}{\end{equation}}

\newcommand{\Bi}{{\cal B}}

\newcommand{\Di}{{\cal D}}

\newcommand{\Hi}{{\cal H}}

\newcommand{\Ki}{{\cal K}}

\newcommand{\Wi}{{\cal W}}

\newcommand{\eps}{\epsilon}

\newcommand{\R}{{\mathbb R}}
\newcommand{\N}{{\mathbb N}}
\newcommand{\Z}{{\mathbb Z}}

\renewcommand{\P}{{\mathbb P}}
\newcommand{\E}{{\mathbb E}}

\newcommand{\Zev}{{\Z^2_{\rm even}}}

\newcommand{\Rc}{R^2_{\rm c}}

\begin{document}

\makeatletter\@addtoreset{equation}{section}
\makeatother\def\theequation{\thesection.\arabic{equation}}

\renewcommand{\labelenumi}{{(\roman{enumi})}}

\title{Brownian Web in the Scaling Limit of Supercritical Oriented Percolation in Dimension $1+1$}
\author{Anish Sarkar$^{\,1}$ \and Rongfeng Sun$^{\,2}$}
\date{Feb 4, 2013}

\maketitle

\footnotetext[1]{Theoretical Statistics and Mathematics Unit, Indian Statistical Institute, New Delhi,
7 S.~J.~S.\ Sansanwal Marg, New Delhi 110016, India. Email: anish@isid.ac.in}

\footnotetext[2]{Department of Mathematics, National University of Singapore, 10 Lower Kent Ridge Road, 119076 Singapore. Email:
matsr@nus.edu.sg}

\begin{abstract}

We prove that, after centering and diffusively rescaling space and time, the collection of rightmost infinite open paths in a supercritical oriented percolation
configuration on the space-time lattice $\Z^2_{\rm even}:=\{(x,i)\in \Z^2: x+i \mbox{ is even}\}$ converges in distribution to the Brownian web. This proves a conjecture of Wu and Zhang~\cite{WZ08}. Our key observation is that each rightmost infinite open path can be approximated by a percolation exploration cluster, and different exploration clusters evolve independently before they intersect.

\end{abstract}

\vspace{12pt}

\section{Introduction}
\subsection{Model and Description of Main Result}\label{S:model}
Let $\Z^2_{\rm even}:=\{(x,i)\in \Z^2: x+i \mbox{ is even}\}$ be a space-time lattice, with oriented edges leading from $(x,i)$ to $(x\pm 1, i+1)$ for all $(x,i)\in \Z^2_{\rm even}$. Oriented percolation on $\Z^2_{\rm even}$ with parameter $p\in [0,1]$ is a random edge configuration on $\Z^2_{\rm even}$, where independently each oriented edge is open with probability $p$, and closed with probability $1-p$. We use $\P_p$ and $\E_p$ to denote respectively probability and expectation for this product probability measure on edge configurations with parameter $p$.

By convention, if there is an open path of oriented edges leading from $z_1=(x_1, i_1)$ to $z_2=(x_2, i_2)$ in $\Z^2_{\rm even}$, then we say that $z_2$ can be reached from $z_1$ and denote it by $z_1 \to z_2$. For any $z\in \Z^2_{\rm even}$, the open cluster at $z$ is then defined by
$$
C_z := \{w\in \Z^2_{\rm even} : z \to w\}.
$$
When $|C_z|$, the cardinality of $C_z$, is infinite, we call $z$ a {\it percolation point}. The set of percolation points will be denoted by $\Ki$.
It is well known (see e.g.~\cite{D84, BG90}) that there exists a critical $p_c\in (0,1)$ such that
$$
\theta(p):=\P_p(|C_{(0,0)}|=\infty)
\ \left\{
\begin{aligned}
=0 & \qquad \mbox{ if } p\in [0, p_c], \\
>0 & \qquad \mbox{ if } p\in (p_c, 1].
\end{aligned}
\right.
$$
The three regimes $p\in [0, p_c)$, $p=p_c$, and $p\in (p_c, 1]$ are called respectively the sub-critcial, critical, and super-critical regimes of oriented percolation.
For a survey of classical results on oriented percolation on $\Z^2_{\rm even}$, see \cite{D84}.

\begin{figure}[tp] 
\begin{center}
\includegraphics[width=9cm]{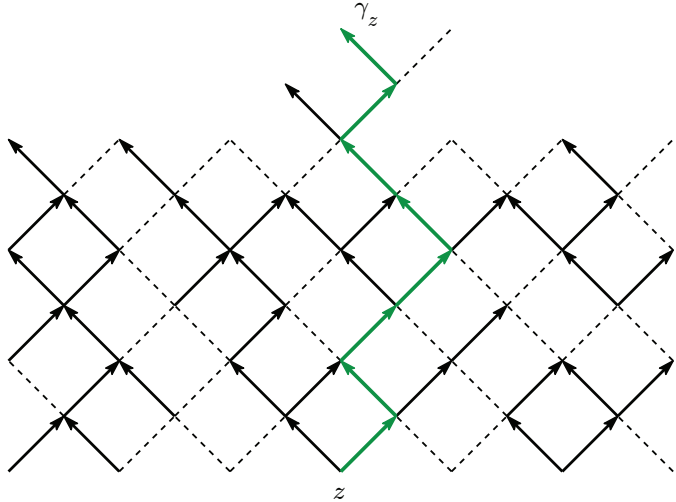}
\caption{An oriented percolation configuration with a rightmost infinite open path $\gamma_z$ drawn in solid green arrows. Closed edges are drawn in dashed lines.}
\label{fig:OPConfig}
\end{center}
\end{figure}

From now on and for the rest of the paper, we restrict our attention to a fixed $p\in (p_c,1)$ and suppress $p$ in $\P_p$ and $\E_p$ to simplify notation. In the supercritical regime, the set of percolation points $\Ki$ is a.s.\ infinite by ergodicity. For each $z=(x,i)\in \Ki$, there is a well-defined {\em rightmost infinite open path} starting from $z$, which we denote by $\gamma_z$ (see Figure~\ref{fig:OPConfig}).
More precisely, $\gamma_z$ can be taken as a mapping from $\{i, i+1, \cdots\}$ to $\Z$ such that $\gamma_z(i)=x$, $(\gamma_z(j), j) \to (\gamma_z(j+1), j+1)$ for all $j\geq i$, and if $\pi$ is any other infinite open path starting from $z$, then $\gamma_z(j)\geq \pi(j)$ for all $j\geq i$. We are interested in the collection of all rightmost infinite open paths in the supercritical oriented percolation configuration:
$$
\Gamma := \{\gamma_z : z\in \Ki\}.
$$

Conditional on $o:=(0,0)\in \Ki$, results of Durrett~\cite{D84} imply that there exists a speed $\alpha:=\alpha(p)>0$ for $p>p_c$, such that
$\lim_{n\to\infty}\gamma_o(n)/n = \alpha$ almost surely. Later, a central limit theorem was established by Kuczek~\cite{K89}, which implies that there exists $\sigma:=\sigma(p)>0$ for $p\in (p_c, 1)$, such that $\frac{\gamma_o(n)-\alpha n}{\sigma \sqrt{n}}$ converges in distribution to a standard normal random variable. As
we will show later, Kuczek's argument further implies that $(\frac{\gamma_o(nt)-\alpha nt}{\sigma \sqrt{n}})_{t\geq 0}$ converges in distribution to a standard Brownian motion. A natural question then arises: If a linear drift $\alpha$ is removed from each path in $\Gamma$, space is rescaled by $\sigma \sqrt{n}$ and time rescaled by $n$,
what is the scaling limit of the whole collection $\Gamma$?

Wu and Zhang~\cite{WZ08} conjectured that the scaling limit of $\Gamma$ should be the so-called {\em Brownian web}, which loosely speaking is a collection of coalescing Brownian motions starting from every point in the space-time plane $\R\times\R$. In~\cite{WZ08}, the authors made the first step towards this conjecture by proving that every pair of paths in $\Gamma$ must coalesce in finite time. Our goal in this paper is to give a proof of this conjecture.

At first sight, it may look surprising that $\Gamma$ (after centering and scaling) should converge to the Brownian web, because of the seemingly complex dependency between paths in $\Gamma$ and the fact that each path depends on the infinite future. However, we make the following key observation which untangles the dependency in a simple way: Each path in $\Gamma$ can be approximated by a {\em percolation exploration cluster} which evolves in a Markovian way, and different exploration clusters evolve independently before they intersect, and ``coalesce'' after they intersect. In the diffusive scaling limit, the width of each cluster tends to 0, while the evolving clusters themselves converge to Brownian motion paths. This justifies heuristically the convergence of $\Gamma$ to the Brownian web. As a byproduct of our approach, we recover the main result in~\cite{WZ08}, that any two paths in $\Gamma$ must coalesce a.s.\ in finite time. We remark that although the heuristic above is simple and natural, some careful analysis is required due to the non-trivial coalescent interaction between exploration clusters after they intersect.

In the remaining subsections of this introduction, we  first recall the characterization of the Brownian web and the relevant topology, and then formulate rigorously our main convergence result. We  then recall the convergence criteria for the Brownian web which we  need to verify, and then end with a discussion on related results and an outline of the rest of the paper.

\subsection{Brownian Web: Characterization}\label{S:topology}

The Brownian web, denoted by $\Wi$, originated from the work of Arratia~\cite{A79, A81} on the scaling limit of the voter model on $\Z$. It arises naturally as the diffusive scaling limit of the dual system of one-dimensional coalescing random walk paths starting from every point on the space-time lattice. We can thus think of the Brownian web as a collection of one-dimensional coalescing Brownian motions starting from every point in the space-time plane $\R^2$, although there is some technical difficulty involved in dealing with an uncountable number of coalescing Brownian motions. Detailed analysis of the Brownian web has been carried out by T\'oth and Werner in~\cite{TW98}. Later, Fontes, Isopi, Newman and Ravishankar~\cite{FINR04} introduced a framework in which the Brownian web is realized as a random variable taking values in the space of {\em compact sets of paths}, which is Polish when equipped with a suitable topology. Under this setup, the object initially proposed by Arratia~\cite{A81} takes on the name the {\em Brownian web}, and we can apply standard theory of weak convergence to prove convergence of various one-dimensional coalescing systems to the Brownian web.

\begin{figure}
\begin{center}
\setlength{\unitlength}{.7cm}
\begin{picture}(10,8)(-5,-4)
\linethickness{.4pt}
\qbezier(0,-3)(0,0)(0,3)
\qbezier(-3,0)(0,0)(3,0)
\linethickness{.4pt}
\qbezier(0,-3)(-6,0)(0,3)
\qbezier(0,-3)(-5.3,0)(0,3)
\qbezier(0,-3)(-4.5,0)(0,3)
\qbezier(0,-3)(-3,0)(0,3)
\qbezier(0,-3)(3,0)(0,3)
\qbezier(0,-3)(4.5,0)(0,3)
\qbezier(0,-3)(5.3,0)(0,3)
\qbezier(0,-3)(6,0)(0,3)
\qbezier(-1.9,-1.8)(0,-1.8)(1.9,-1.8)
\qbezier(-2.2,-1.5)(0,-1.5)(2.2,-1.5)
\qbezier(-2.65,-1)(0,-1)(2.65,-1)
\qbezier(-2.65,1)(0,1)(2.65,1)
\qbezier(-2.2,1.5)(0,1.5)(2.2,1.5)
\qbezier(-1.9,1.8)(0,1.8)(1.9,1.8)

\put(0,-3){\circle*{.25}}
\put(0,3){\circle*{.25}}
\put(0,0){\circle*{.25}}
\put(2.2,1.5){\circle*{.25}}
\put(-2.65,-1){\circle*{.25}}

\put(-.6,-3.5){$(\ast,-\infty)$}
\put(-.6,3.5){$(\ast,+\infty)$}
\put(-1.3,-.5){$(0,0)$}
\put(2.3,1.6){$(+\infty,2)$}
\put(-5.3,-1.1){$(-\infty,-1)$}

\end{picture}
\caption[The compactification $\Rc$ of $\R^2$.]{The compactification
$\Rc$ of $\R^2$.}\label{fig:comp}
\end{center}
\end{figure}
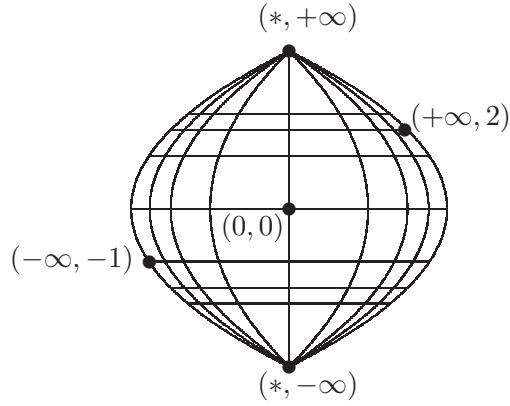

We now recall from~\cite{FINR04} the space of {\em compact sets of paths} in which the Brownian web $\Wi$ takes its value. Let $\Rc$ denote the completion of the space-time plane $\R^2$ w.r.t.\ the metric
\be\label{rho}
\rho\big((x_1, t_1), (x_2,t_2)\big) = \left|\tanh(t_1)-\tanh(t_2)\right|
\ \vee\ \left|\frac{\tanh(x_1)}{1+|t_1|}-\frac{\tanh(x_2)}{1+|t_2|}\right|.
\ee
As a topological space, $\Rc$ can be identified with the continuous image of
$[-\infty, \infty]^2$ under a map that identifies the line
$[-\infty,\infty]\times\{\infty\}$ with a single point $(*, \infty)$, and the
line $[-\infty,\infty]\times\{-\infty\}$ with the point $(*,-\infty)$, see
Figure~\ref{fig:comp}.

A path $\pi$ in $\Rc$, whose starting time we denote by
$\sigma_\pi\in [-\infty,\infty]$, is a mapping $\pi :
[\sigma_\pi,\infty] \to [-\infty, \infty] \cup\{*\}$ such that
$\pi(\infty)=*$, $\pi(\sigma_\pi)=*$ if $\sigma_\pi=-\infty$, and $t
\to (\pi(t), t)$ is a continuous map from $[\sigma_\pi,\infty]$ to
$(\Rc, \rho)$. We then define $\Pi$ to be the space of all
paths in $\Rc$ with all possible starting times in $[-\infty,\infty]$.
Endowed with the metric
\be\label{PId}
\!
d(\pi_1, \pi_2) = \Big|\!\tanh(\sigma_{\pi_1})-\tanh(\sigma_{\pi_2})\Big|
\ \vee \sup_{t\geq \sigma_{\pi_1} \wedge \sigma_{\pi_2}}\!
\left|\frac{\tanh(\pi_1(t\vee \sigma_{\pi_1}))}{1+|t|}
-\frac{\tanh(\pi_2(t\vee \sigma_{\pi_2}))}{1+|t|}\right|, \!\!\!
\ee
$(\Pi, d)$ is a complete separable metric space. Note that
convergence in the metric $d$ can be desrcibed as locally uniform convergence
of paths plus convergence of starting times. (The metric
$d$ differs slightly from the original choice in \cite{FINR04}, which
is somewhat less natural as explained in the appendix of \cite{SS08}.)

We can now define $\Hi$, the {\em space of compact subsets of $(\Pi, d)$},
equipped with the Hausdorff metric
\be\label{dH}
d_{\Hi}(K_1, K_2) = \sup_{\pi_1\in K_1} \inf_{\pi_2\in K_2}\!\!
d(\pi_1, \pi_2)\ \vee  \sup_{\pi_2\in K_2}
\inf_{\pi_1\in K_1} d(\pi_1, \pi_2).
\ee
The space $(\Hi, d_\Hi)$ is also a complete separable metric
space. Let $\Bi_\Hi$ be the Borel $\sigma$-algebra associated with
$d_\Hi$. The Brownian web $\Wi$ is  an
$(\Hi, \Bi_\Hi)$-valued random variable.

Following convention, for $K\in\Hi$ and $A\subset\Rc$, we  let
$K(A)$ denote the set of paths in $K$ with starting points in
$A$. When $A=\{z\}$ for $z\in\Rc$, we also write $K(z)$ instead
of $K(\{z\})$.

We now recall from \cite[Theorem 2.1]{FINR04} the following characterization of the
Brownian web $\Wi$.
\bt\label{T:webchar}{\bf[Characterization of the Brownian web]}
There exists an $(\Hi, \Bi_\Hi)$-valued random variable $\Wi$, called
the standard Brownian web, whose distribution is uniquely determined
by the following properties:
\begin{itemize}
\item[{\rm (a)}] For each deterministic $z\in\R^2$, almost surely
there is a unique path $\pi_z\in \Wi(z)$.

\item[{\rm (b)}] For any finite deterministic set of points $z_1, \ldots, z_k
  \in\R^2$, the collection $(\pi_{z_1}, \ldots, \pi_{z_k})$
  is distributed as coalescing Brownian motions.

\item[{\rm (c)}] For any deterministic countable dense subset
$\Di \subset\R^2$, almost surely, $\Wi$ is the closure of
$\{\pi_z : z\in\Di\}$ in $(\Pi, d)$.
\end{itemize}
\et
Theorem~\ref{T:webchar} shows that the Brownian web is in some sense {\em separable}: even
though there are uncountably many coalescing Brownian motions in the Brownian web, the
whole collection is a.s.\ determined uniquely by a countable skeletal subset of paths.

\subsection{Formulation of Main Result}
We formulate in this subsection the convergence of $\Gamma$, the collection of rightmost
infinite open paths, to the Brownian web $\Wi$ after suitable centering and scaling.

Given a fixed $p\in (p_c, 1)$, let $\alpha:=\alpha(p)>0$
and $\sigma:=\sigma(p)>0$ be as introduced in Section~\ref{S:model}, such that conditional on $o:=(0,0)$
being a percolation point, $\frac{\gamma_o(n)-\alpha n}{\sigma \sqrt{n}}$ converges in distribution to a standard
normal. We will formulate this convergence precisely in Lemma~\ref{L:KCLT}, where we recall Kuczek's proof of the
central limit theorem and extend it for our purposes.

For each percolation point $z=(x,i)\in \Ki$, we first extend the definition of the rightmost infinite open
path $\gamma_z$ from the domain $\{i, i+1, \ldots\}$ to $[i, \infty]$ such that $\gamma_z$ interpolates linearly
between consecutive integer times and $\gamma_z(\infty)=*$. With this extended definition of $\gamma_z$, which
we still denote by $\gamma_z$ for convenience, it becomes a path in the space $(\Pi, d)$ introduced in
Section~\ref{S:topology}. We will then let $\Gamma:=\{\gamma_z: z\in\Ki\}$ denote the set of extended rightmost
infinite open paths in the percolation configuration. Since paths in $\Gamma$ are a.s.\ equicontinuous, $\overline{\Gamma}$,
the closure of $\Gamma$ in $(\Pi, d)$, is a.s.\ compact and hence $\overline{\Gamma}$ is a random variable
taking values in $(\Hi, \Bi_{\Hi})$, the space of compact subsets of $(\Pi,d)$. Note that $\overline{\Gamma}\backslash\Gamma$
only contains paths of the form $\pi: [\sigma_\pi, \infty]\to [-\infty, \infty]\cup\{*\}$ with either $\sigma_\pi \in\R$ and
$\pi(t)\equiv \pm \infty$ for all $t\geq \sigma_\pi$; or $\sigma_\pi=\infty$; or $\sigma_\pi=-\infty$, in which case for any
$t>-\infty$, there exists some $\gamma\in \Gamma$ such that $\pi=\gamma$ on $[t,\infty]$. In other words, taking
the closure of $\Gamma$ in $(\Pi, d)$ does not alter the configuration of paths in $\Gamma$ restricted to any finite space-time
region. Therefore it suffices to study properties of $\Gamma$ instead of $\overline\Gamma$ in our analysis.

To remove a common drift from all paths in $\Gamma$ and perform diffusive scaling of space and time,  we define for any
$a\in\R$, $b, \eps>0$, a shearing and scaling map $S_{a,b,\eps}: \Rc\to\Rc$ with
\be\label{Smap}
S_{a,b,\eps} (x,t) :=
\begin{cases}
\big(\frac{\sqrt{\eps}}{b}(x-at), \eps t\big) & \mbox{if }  (x,t) \in \R^2, \\
(\pm \infty, \eps t) & \mbox{if } (x,t)=(\pm \infty, t) \mbox{ with } t\in\R, \\
(*, \pm \infty)  & \mbox{if } (x,t)=(*, \pm \infty),
\end{cases}
\ee
where $a$ is the drift that is being removed by a shearing of $\Rc$, $\eps$ is the diffusive scaling parameter, and $b$ determines
the diffusion coefficient in the diffusive scaling. When $t$ is understood to be a time, we will define
\be\label{Smapt}
S_{a,b,\eps}t:=\eps t.
\ee
Note that $S_{a,b,\eps}$ can be obtained by first applying the shearing map
$S_{a, 1,1}$ and then the diffusive scaling map $S_{0, b, \eps}$. By identifying a path $\pi \in \Pi$ with its graph in $\Rc$, we can also define
$S_{a,b,\eps}: (\Pi, d)\to (\Pi, d)$ by applying $S_{a,b,\eps}$ to each point on the graph of $\pi$. Similarly, if $K\subset \Pi$,
then $S_{a,b,\eps}K := \{S_{a,b,\eps}\pi : \pi \in K\}$. If $K\in \Hi$, then it is clear that also $S_{a,b,\eps}K \in \Hi$. Therefore
$S_{\alpha, \sigma, \eps}\overline{\Gamma}$ is also an $(\Hi, \Bi_{\Hi})$-valued random variable.

We can now formulate the main result of this paper.
\bt\label{T:main}{\bf [Convergence to the Brownian web]} Let $p\in (p_c, 1)$ and let $\overline\Gamma$ be defined as above. There exist
$\alpha, \sigma>0$ such that as $\eps\downarrow 0$, the sequence of $(\Hi, \Bi_{\Hi})$-valued random variables
$S_{\alpha, \sigma, \eps}\overline{\Gamma}$ converges in distribution to the standard Brownian web $\Wi$.
\et

\subsection{Brownian Web: Convergence Criteria}
We will prove Theorem~\ref{T:main} by verifying the convergence criteria for the Brownian web proposed in~\cite{FINR04},
which we now recall.

For a compact set of paths $K\in \Hi$, and for $t>0$ and $t_0,a,b\in\R$ with $a<b$, let
\be\label{eta}
\eta_K(t_0,t;a,b) := \big|\{\pi(t_0+t)\, :\, \pi \in K \mbox{ with } \sigma_\pi\leq t_0 \mbox{ and } \pi(t_0)\in [a,b]  \}\big|,
\ee
which counts the number of distinct points on $\R\times\{t_0+t\}$ touched by some path in $K$ which also touches $[a,b]\times\{t_0\}$.

An $(\Hi, \Bi_{\Hi})$-valued random variable ${\cal X}$ is said to have {\em non-crossing paths} if a.s.\ there exist no $\pi, \tilde\pi \in {\cal X}$
such that $(\pi(t)-\tilde\pi(t))(\pi(s)-\tilde\pi(s))<0$ for some $s,t\geq \sigma_\pi \vee \sigma_{\tilde\pi}$. Note that $\overline{\Gamma}$ has non-crossing path. For $(\Hi, \Bi_{\Hi})$-valued random variables with non-crossing paths, the following convergence criteria
was formulated in~\cite[Theorem 2.2]{FINR04}.

\bt\label{T:convct1}{\bf [Convergence criteria]} Let $({\cal X}_n)_{n\in\N}$ be a sequence of $(\Hi, \Bi_{\Hi})$-valued random variables with
non-crossing paths. If the following conditions are satisfied, then ${\cal X}_n$ converges in distribution to the standard Brownian web $\Wi$.
\begin{itemize}
\item[{\rm \,(I)\ }] Let $\Di$ be a deterministic countable dense subset of $\R^2$. Then there exist $\pi^y_n \in {\cal X}_n$ for $y\in\Di$ such that, for each
finite collection $y_1, y_2, \ldots, y_k\in \Di$, $(\pi^{y_1}_n, \ldots, \pi^{y_k}_n)$ converge in distribution as $n\to\infty$ to a collection of coalescing Brownian motions starting at $(y_1,\ldots, y_k)$.

\item[{\rm (B1)}] For all $t>0$, $\limsup_{n\to\infty} \sup_{(a,t_0)\in\R^2} \P(\eta_{{\cal X}_n}(t_0, t; a, a+\eps)\geq 2) \to 0$ as $\eps\downarrow 0$.
\item[{\rm (B2)}] For all $t>0$, $\eps^{-1}\limsup_{n\to\infty} \sup_{(a,t_0)\in\R^2} \P(\eta_{{\cal X}_n}(t_0, t; a, a+\eps)\geq 3) \to 0$ as $\eps\downarrow 0$.
\end{itemize}
\et
As shown in \cite[Prop.~B.2]{FINR04}, condition (I) and the non-crossing property imply that $({\cal X}_n)_{n\in\N}$ is a tight sequence of $(\Hi, \Bi_{\Hi})$-valued
random variables. Condition (I) also guarantees that any subsequential weak limit of $({\cal X}_n)_{n\in\N}$ contains as many paths as, possibly more than,
the Brownian web $\Wi$. Conditions (B1) and (B2) are density bounds which rule out the presence of extra paths other than the Brownian web paths in any
subsequential weak limit.

As alluded to at the end of Section~\ref{S:model}, we will verify condition (I) by approximating each path in $\Gamma$ by a percolation exploration cluster which enjoys Markov and independence properties. The verification of (B1) is closely related to that of (I). The verification of condition (B2) typically relies on FKG inequalities for the law of each individual path in ${\cal X}_n$ (see e.g.~\cite[Theorem 6.1]{FINR04}, \cite[Lemma 2.6]{FFW05}, and \cite[Section 2.1]{CFD09}). Although we will be replacing each path in $\Gamma$ by an exploration cluster, it turns out that we can still apply FKG for the underlying percolation edge configuration to deduce (B2). We remark that there is an alternative convergence criterion formulated in \cite[Theorem 1.4]{NRS05}, which is often easier to verify than (B2) when FKG inequalities are not applicable. In fact we first proved Theorem~\ref{T:main} by verifying the convergence criteria in \cite[Theorem 1.4]{NRS05} without using FKG (see~\cite{SS12}). The argument is lengthier and more involved, but in a sense more robust.

\subsection{Discussion and Outline}
The Brownian web arises as the diffusive scaling limit of many one-dimensional coalescing systems. The prime example is the collection of coalescing simple random walks
on $\Z$, for which the convergence to the Brownian web was established in~\cite[Theorem 6.1]{FINR04}. This result was extended to general coalescing random walks
with crossing paths in~\cite{NRS05} under a finite $5$-th moment assumption on the random walk increment, which was later improved in~\cite{BMSV06} to an essentially optimal assumption of finite $(3+\eps)$-th moment for any $\eps>0$.

Other one-dimensional coalescing systems (all with non-crossing paths) which have been shown to converge to the Brownian web include: two-dimensional Poisson trees, which was introduced in~\cite{FLT04} and shown to converge to the Brownian web in~\cite{FFW05}; a two-dimensional drainage network model which was introduced in~\cite{GRS04} and shown to converge to the Brownian web in~\cite{CFD09}, as well as an extension studied more recently in~\cite{CV11}. Interestingly, the Brownian web, or rather, the coalescing flow generated by it known as the {\em Arratia flow}, also arises in the scaling limit of a planar aggregation model, see~\cite{NT11a, NT11b}, where convergence was established using a different topology tailored more specifically for the study of stochastic flows.

Another one-dimensional coalescing system conjectured to converge to the Brownian web is the {\em directed spanning forest}, which was introduced in~\cite{BB07} and shown recently to be a.s.\ a tree in~\cite{CT11}. As one might expect, the difficulty in establishing convergence to the Brownian web
lies in the specific form of dependence in that model.

Instead of considering the collection of rightmost infinite open paths in a supercritical oriented percolation configuration, one may also fix a realization of the percolation configuration and consider the set of directed coalescing random walk paths on the infinite clusters. Namely, from each percolation point $z=(x,t)\in\Ki$, a random walk starts from $z$ and jumps to a site in $\{(x+1, t+1), (x-1, t+1)\}\cap \Ki$ with uniform probability, and different random walks coalesce when they meet. Naturally one expects such a collection of coalescing random walk paths to converge to the Brownian web under diffusive scaling, both for the quenched measure of the random walks under a typical realization of the percolation configuration, and for the averaged (or annealed) measure where the law of the percolation configuration is integrated out. The percolation exploration procedure we devise in this paper can also be used to study this model. Kuczek's proof of the central limit theorem~\cite{K89} for a single path in $\Gamma$, say $\gamma_o$, is based on the identification of so-called {\em break points} along $\gamma_0$, which gives a renewal decomposition of $\gamma_0$. Our exploration procedure (which always explores the two outgoing edges from each $z=(x,t)\in\Zev$ from right to left) provides a systematic way of discovering these break points. If we instead explore the two outgoing edges from each $z\in\Zev$ with a random order, then we obtain exactly a directed random walk on the supercritical oriented percolation clusters, with a renewal decomposition of the random walk path under the averaged measure. The exploration procedure ensures that different random walk paths evolve independently (under the averaged measure) before their associated exploration clusters intersect. We learned recently from Matthias Birkner that the extension of Kuczek's idea of break points to a directed random walk on supercritical oriented percolation clusters in fact dates back to Neuhauser~\cite{N92}. This idea was used recently by Birkner et al~\cite{BCDG12} to prove averaged and quenched invariance principles for a directed random walk on supercritical oriented percolation clusters.

Our result sheds some light on the scaling limit of super-critical oriented percolation in dimension 1+1. We expect the same result to hold for the one-dimensional
contact process. However, a most interesting and challenging direction of extension will be to investigate what kind of structures appear in the scaling limit of {\em critical} oriented percolation, or {\em near-critical} oriented percolation where $p\downarrow p_c$ in tandem with the scaling of
the lattice. For two-dimensional unoriented percolation, the critical scaling limit has been constructed and studied using Schramm-L\"owner Evolutions, while the
near-critical scaling limit is currently under construction (see~\cite{GPS10} and the references therein). However for oriented percolation, there is still no conjecture on how to characterize the critical and near-critical scaling limits, or if such limits exist at all.

The rest of the paper is organized as follows. In Section~\ref{S:explore}, we define and establish some basic properties for the exploration clusters which approximate the rightmost infinite open paths. We will also recall Kuczek's proof~\cite{K89} of the central limit theorem for a rightmost infinite  open path and extend it
to establish an invariance principle. In Section~\ref{S:I}, we will prove the convergence of multiple exploration clusters to coalescing Brownian motions, which implies condition (I). Lastly in Section~\ref{S:B12}, we will verify conditions (B1) and (B2), which completes the proof of Theorem~\ref{T:main}.

\section{Exploration Clusters}\label{S:explore}
In this section, we introduce the key objects in our analysis, the {\em percolation exploration clusters}. We will establish some basic properties for these exploration clusters. We  then show that each path in $\Gamma$ can be approximated by an associated exploration cluster, in the sense that both converge to the same Brownian motion after suitable centering and scaling.

\subsection{Construction of exploration clusters}
So far the rightmost infinite open paths $\gamma_z$ are only defined for $z\in \Ki$,
the set of percolation points. We first extend this definition to all $z=(x,i)\in \Z^2_{\rm even}$ by defining
\be\label{gammaz}
\gamma_z:=\gamma_{z'}, \quad \mbox{where } z'=(y,i) \mbox{ with } y=\max\{u\leq x : (u,i) \in \Ki\},
\ee
which is the rightmost infinite open path starting from $(-\infty, x]\times\{i\}$, and it exists a.s.\ because we assumed that $p>p_c$. Note that $\Gamma=\{\gamma_z: z\in \Ki\}=\{\gamma_z: z\in\Z^2_{\rm even}\}$.

Without loss of generality, let $z=o$ be the origin. In a nutshell, the exploration cluster $C_o(n)$ we use to approximate $\gamma_o$ up to time $n$, consists of the minimal set of open and closed edges, denoted respectively by $E^{\rm o}_o(n)$ and $E^{\rm c}_o(n)$, that need to be explored in order to find the {\em rightmost open path} connecting $(-\infty, 0]\times\{0\}$ to $\Z\times\{n\}$.

\begin{figure}[tp] 
\begin{center}
\includegraphics[width=13cm]{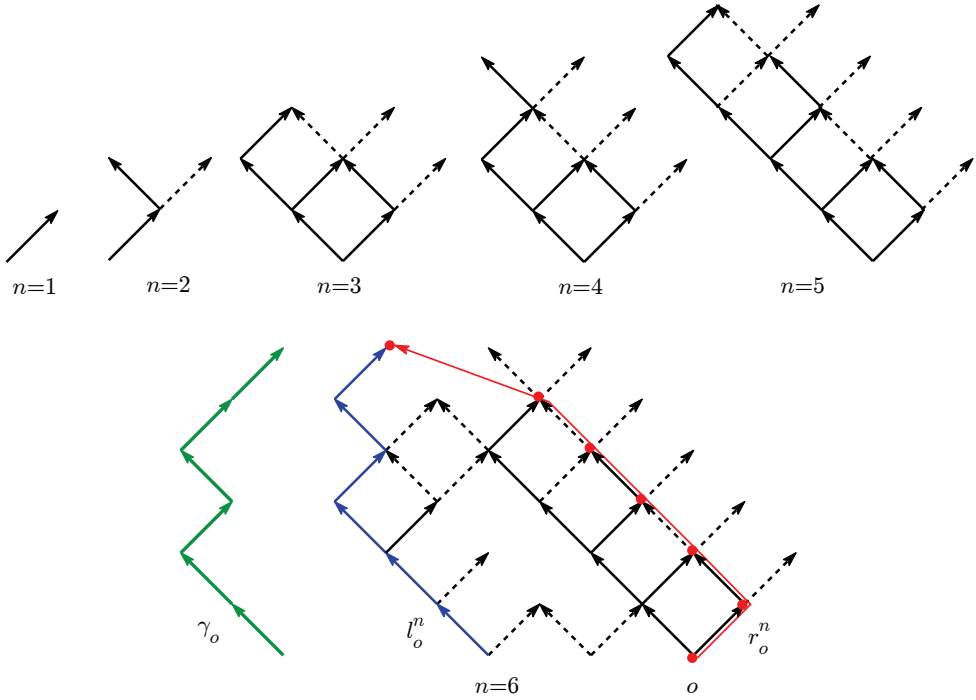}
\caption{The exploration cluster process $C_o(n)$ for $n=1,2,\ldots, 6$.}
\label{fig:exploration}
\end{center}
\end{figure}

To construct $C_o(n):=(E^{\rm o}_o(n), E^{\rm c}_o(n))$, we start by exploring the open cluster at the origin in the following way (see Figure~\ref{fig:exploration}). 
Let $e_o^+$ and $e_o^-$ denote respectively the up-right and up-left edge starting from $o$, and let $w^+$ and $w^-$ denote respectively the up-right and up-left neighbor of $o$. We first explore the status of $e_o^+$. If $e_o^+$ is closed, then we move on to explore $e_o^-$. If $e_o^+$ is open, then we explore next the open cluster at $w^+$ by the same procedure and stop when either the first open path reaching $\Z\times\{n\}$ is discovered, or the exploration of the open cluster at $w^+$ is finished without discovering an open path to $\Z\times\{n\}$, in which case we move on to explore $e_o^-$. When it comes to exploring $e_o^-$, if it is closed, then the exploration of the open cluster at $o$ is finished; otherwise we explore next the open cluster at $w^-$ by the same procedure until either the first open path reaching $\Z\times\{n\}$ is discovered, or the exploration of the open cluster at $w^-$ is finished without discovering an open path to $\Z\times\{n\}$, in which case the exploration of the open cluster at $o$ is again finished. We explore the open clusters at $o=(0,0)$, $(-2,0)$, $(-4,0)$, $\ldots$ in this order, until a first open path connecting $(-\infty, 0]\times\{0\}$ and $\Z\times\{n\}$ is discovered, at which time we stop the exploration. The sets $E^{\rm o}_o(n)$ and $E^{\rm c}_o(n)$ are the sets of open and respectively closed edges that have been explored up to the time the exploration process is stopped.

By construction, there is a unique open path using edges in $E^{\rm o}_o(n)$ which connects $(-\infty, 0]\times\{0\}$ to $\Z\times\{n\}$,
and it is also the rightmost open path among all open paths connecting $(-\infty, 0]\times\{0\}$ to $\Z\times\{n\}$. We denote this path by $l_o^{n} :\{0,1,\ldots, n\} \to \Z$. Let $V_o(n)$ denote the set of vertices at which at least one of the two outgoing edges have been explored. The path $l_o^n$ serves as both the left vertex boundary and the left edge boundary of $C_o(n)$ in the senese that all vertices in $V_o(n)$ and all edges in $E_o^{\rm o}(n)\cup E_o^{\rm c}(n)$ must lie on or to the right of $l^n_o$. We can similarly define the right vertex boundary of $C_o(n)$ by the path
\be
r_o^n(j) := \max\{x\in \Z: (-\infty, 0]\times\{0\} \to (x, j)\}, \qquad j=0,1,\ldots, n.
\ee
Note that $l^n_o(n)=r_o^n(n)$ for all $n\geq 0$. We will call $(l_o^n, r_o^n)$ the {\em left and right boundaries} associated with the exploration cluster $C_o(n)$. The time-evolution of $(l_o^n, r_o^n)_{n\geq 0}$ will be called the {\em left and right boundary processes} associated with the {\em exploration cluster process} $C_o:=(C_o(n))_{n\geq 0}$. For a general $z=(x,i)\in\Z^2_{\rm even}$, the exploration cluster process $C_z:=(C_z(n))_{n\geq i}$ and its left and right boundary processes $(l_z^n, r_z^n)_{n\geq i}$ are defined similarly.

Note that the exploration process in the construction of $C_o(n)$, $n\geq 0$, discovers the percolation configuration on $\Z^2_{\rm even}$ in a Markovian way: conditional on the status of the edges that have been explored, the next edge to be explored is determined uniquely, and it is conditionally open with probability $p$ and closed with probability $1-p$. In particular, $(C_o(n))_{n\geq 0}$ is Markov. Furthermore, different exploration cluster processes evolve independently as long as the sets of explored edges do not intersect.

\subsection{Approximation by an exploration cluster}
We now show that for each $z=(x,i)\in\Z^2_{\rm even}$, $\gamma_z\in\Gamma$ can be approximated on the diffusive scale by the exploration cluster process $C_z:=(C_z(n))_{n\geq i}$, or rather, by the associated boundary processes $(l^n_z, r^n_z)_{n\geq i}$. First we collect some basic properties of $\gamma_z$ and $(l^n_z, r^n_z)_{n\geq i}$. Without loss of generality, assume $z=o$. All the discussions and results that will follow and the notation we will introduce adapt straightforwardly to a general $z\in\Z^2_{\rm even}$.

First, we identify $r_o^n: \{0, 1,\ldots, n\}\to \Z$ with its extended definition on $[0,\infty]$ by setting $r_o^n(s):=r_o^n(n)$ for all $s\in [n,\infty)$, $r_o^n(\infty)=*$, and linearly interpolating between consecutive integer times. The same applies to $l_o^n$.

In the construction of the exploration cluster process $C_o:=(C_o(n))_{n\geq 0}$, we observed that $l_o^n$ is the rightmost open path connecting $(-\infty, 0]\times \{0\}$ to $\Z\times\{n\}$, and $r_o^n(j)$, for each $0\leq j\leq n$, is the rightmost position at time $j$ that can be reached  by any open path starting from $(-\infty, 0]\times \{0\}$, while $\gamma_o$ is the rightmost infinite open path starting from $(-\infty, 0]\times\{0\}$. These facts readily imply

\bl\label{L:lrprop} Let $\gamma_o$ be defined as in (\ref{gammaz}), and let $C_o$ be the associated exploration cluster process with left and right boundaries $(l_o^n, r_o^n)_{n\geq 0}$. Then
\begin{itemize}
\item[\rm (i)] There exists $r_o: [0,\infty]\to \Z\cup\{*\}$ such that for all $n\geq 0$, $r_o^n(\cdot) = r_o(\cdot)$ on $[0,n]$, and
$$
r_o(n)= \max\{y\in \Z: (-\infty, 0]\times\{0\} \to (y,n) \}.
$$

\item[\rm (ii)] For all $n\geq 0$, we have $\gamma_o(\cdot) \leq l_o^n(\cdot)\leq r_o(\cdot)$ on $[0,n]$.

\item[\rm (iii)] For all $n\geq 0$, we have $l_o^n(n) = r_o^n(n)$.

\item[\rm (iv)] For all $m\geq n\geq 0$, we have $l_o^m(\cdot)\leq l_o^n(\cdot)$ on $[0, n]$.
\end{itemize}
\el
The time-consistency of $(r^n_o)_{n\geq 0}$ established in Lemma~\ref{L:lrprop}~(i) allows us to replace $(r^n_o)_{n\geq 0}$ by a single path $r_o$.
The left-boundary process $(l^n_o)_{n\geq 0}$ does not share this time-consistency property, as illustrated in Figure~\ref{fig:exploration}. We also note that $\gamma_o$ and $(l_o^n)_{n\geq 0}$ are nearest-neighbor paths, while $r_o$ may have jumps of size more than 1 to the left, but jumps to the right are always nearest-neighbor.
\medskip

By the ordering relation in Lemma~\ref{L:lrprop}~(ii), to show that $\gamma_o$ can be approximated by $(C_o(n))_{n\geq 0}$, which converges to a Brownian motion after
proper centering and scaling, it suffices to show that
\bprop\label{P:Donsker1}{\bf [Invariance Principle]} Let $p\in (p_c, 1)$. There exist $\alpha:=\alpha(p)>0$ and $\sigma:=\sigma(p)>0$, such that as $\eps\downarrow 0$, $S_{\alpha, \sigma, \eps}(\gamma_o, r_o)$ converges in distribution as a sequence of $\Pi\times\Pi$-valued random variables to $(B, B)$ for a standard Brownian motion $B:=(B_t)_{t\geq 0}$, where $S_{\alpha, \sigma, \eps}$ is the shearing and diffusive scaling map defined in (\ref{Smap}).
\eprop

In~\cite[Section 3]{D84}, Durrett considered the very same process $(r_o(n))_{n\geq 0}$ and used the sub-additive ergodic theorem to show that a.s.\
$$
\lim_{n\to\infty} \frac{r_o(n)}{n} =: \alpha =\inf_{n\geq 1} \E\Big[\frac{r_o(n)}{n}\Big] >0 \qquad \mbox{if } p>p_c.
$$
Kuczek~\cite{K89} later established a central limit theorem for a variant of $r_o(n)$. When the cluster at the origin dies out at time $n$, instead of exploring next the cluster at $(-2,0)$ as in our construction of $r_o$, Kuczek explores next the cluster at $(n,n)$ and iterates this process.
Since we are interested in an invariance principle for $r_o$ and $\gamma_o$, as well as a bound on the difference between $r_o$ and $\gamma_o$ which we will need later, we recall below Kuczek's argument and adapt it to prove Proposition~\ref{P:Donsker1}.

A key tool in Kuczek's argument is the use of {\em break points}, which are analogous to regeneration times in the study of random walks in random environments. For $r_o$, they are the successive percolation points in the sequence $(r_o(n), n)_{n\geq 0}$, which we denote by $(r_o(T_i), T_i)_{i\in\N}$. The break points are exactly the points at which $\gamma_o$ and $r_o$ coincide, and it is easy to see that for each $i\in\N$, $l^n_o(\cdot)=\gamma_o(\cdot)$ on $[0, T_i]$ for all $n\geq T_i$. Let
\be\label{Xtau}
\begin{aligned}
\tau_1& :=T_1,   & \tau_i& :=T_i-T_{i-1}  &\mbox{ for } i\geq 2; \\
X_1 & :=r_o(T_1),  & X_i & :=r_o(T_i)-r_o(T_{i-1})  &\mbox{ for } i\geq 2.
\end{aligned}
\ee
Note that when $o$ is a percolation point, $X_1=\tau_1=0$. Kuczek proved in~\cite{K89} that
\bl{\bf [Conditional CLT]} \label{L:KCLT} Conditional on the event $o\in \Ki$, $(X_i, \tau_i)_{i\geq 2}$ are i.i.d.\ with all moments finite, and $(r_o(n)-\alpha n)/\sigma \sqrt{n}$ converges in distribution to the standard Normal random variable as $n\to\infty$, where\footnote{There was a typo in \cite{K89} after Lemma 2, and in \cite[Prop.~2.1]{WZ08}, where  the factor $\E[\tau_2]^{-3}$ was missing in the formulae for $\sigma^2$. }
\be\label{alphasigma}
\alpha =\frac{\E[X_2]}{\E[\tau_2]} \quad \mbox{and} \quad
\sigma^2 =\frac{\E\big[(X_2\E[\tau_2]-\tau_2\E[X_2])^2\big]}{\E[\tau_2]^3}>0.
\ee
\el
 Kuczek's proof of Lemma~\ref{L:KCLT} is based
on the key observation that conditional on $z=(x,i)$ being a break point along $r_o$, the percolation configurations before and after time $i$ are independent.

Kuczek's arguments can be extended to prove an invariance principle.
\bl{\bf [Conditional Invariance Principle]}\label{L:KDonsker} Conditional on the event $o\in \Ki$, $S_{\alpha, \sigma, \eps}r_o$ converges in distribution as a sequence of $\Pi$-valued random variables to
a standard Brownian motion $B:=(B_t)_{t\geq 0}$, where $\alpha$ and $\sigma$ are as in (\ref{alphasigma}).
\el
{\bf Proof.} First we replace $(r_o(t))_{t\geq 0}$ by a path $(\tilde r_o(t))_{t\geq 0}$, where $\tilde r_o(T_i):=r_o(T_i)$ for all $i\in\N$,
and for $t\in (T_i, T_{i+1})$, $\tilde r_o(t)$ is defined by linearly interpolating between $r_o(T_i)$ and $r_o(T_{i+1})$. Then
$$
S_{\alpha, \sigma, \eps}r_o = S_{\alpha, \sigma, \eps}\tilde r_o + S_{0,\sigma,\eps}(r_o-\tilde r_o),
$$
and it suffices to show that $S_{0,\sigma,\eps}(r_o-\tilde r_o)$ converges in distribution to the zero function, while $S_{\alpha, \sigma, \eps}\tilde r_o$ converges in
distribution to $B$.

If we let $\vec X_i := (r_o(T_{i-1}+j)-r_o(T_{i-1}))_{0\leq j\leq \tau_i} \in \cup_{n\geq 1} \Z^n$
for $i\geq 2$, then Kuczek's observation on break points also implies that $(\vec X_i, \tau_i)_{i\geq 2}$ are i.i.d.\ conditional on $o\in \Ki$.
In particular,
$$
\sup_{t\in [T_{i-1}, T_i]}|r_o(t)-\tilde r_o(t)|, \quad i\geq 2,
$$
are also i.i.d.\ conditional on $o\in \Ki$. Furthermore, the facts that $r_o(T_{i-1})=\gamma_o(T_{i-1})$ for each $i\geq 2$, $\gamma_o\leq r_o$, $\gamma_o$ is a nearest-neighbor path, while jumps of $r_o$ to the right are always nearest-neighbor, together imply that
$$
\sup_{t\in [T_{i-1}, T_i]}|r_o(t)-\tilde r_o(t)| \leq 2 (T_i-T_{i-1})=2\tau_i, \quad i\geq 2,
$$
which has all moments finite. It is then an easy exercise to verify, which we leave to the reader, that $S_{0,\sigma, \eps}(r_o-\tilde r_o)$
converges in distribution to the zero function as $\eps\downarrow 0$.

Since $S_{\alpha, \sigma, \eps} \tilde r_0= S_{0, \sigma, \eps}S_{\alpha, 0, 0}\tilde r_o$, we first apply the shearing map and note that
$S_{\alpha, 0, 0}\tilde r_o$ is the path obtained by linearly interpolating between the sequence of space-time points $(\sum_{i=1}^n \tilde X_i, \sum_{i=1}^n \tau_i)$, $n\geq 1$, where $\tilde X_i=X_i-\alpha \tau_i$. We can therefore regard $S_{\alpha, 0, 0}\tilde r_o$ as a time change of the random walk $W(n)=\sum_{i=1}^n \tilde X_i$,
with the time change given by $T(n):=\sum_{i=1}^n \tau_i$, and $W(t)$ and $T(t)$ for non-integer $t$ are defined by linearly interpolating between consecutive integer times. More precisely, $(S_{\alpha, 0, 0}\tilde r_o)(t) = W(T^{-1}(t))$. Note that $(\tilde X_i)_{i\geq 2}$ are i.i.d.\ with $\E[\tilde X_i]=0$ and $\E[\tilde X_i^2]=\E[\tau_i] \sigma^2$. Therefore $S_{0, \sigma, \eps}W$ converges in distribution to $\sqrt{\E[\tau_2]}B$ for a standard Brownian motion $B$. On the other hand,
we note that $(\eps T(t/\eps))_{t\geq 0}$ satisfies a law of large numbers and converges in distribution to the linear function $g(t)=\E[\tau_2]t$, $t\geq 0$, with the
topology of local uniform convergence. Therefore $(\eps T^{-1}(t/\eps))_{t\geq 0}$ converges in distribution to $g^{-1}(t)=\E[\tau_2]^{-1} t$, $t\geq 0$. It is then a standard exercise, which we again leave to the reader, to show that
$$
S_{0,\sigma, \eps} S_{\alpha, 0, 0}\tilde r_o = S_{0,\sigma, \eps} W(T^{-1}(t)) = \Big( \frac{\sqrt{\eps} W\big(\eps^{-1} \ \eps T^{-1}(t/\eps)\big)}{\sigma} \Big)_{t\geq 0}
$$
converges weakly to $\sqrt{\E[\tau_2]} B(t/\E[\tau_2]) \stackrel{d}{=} B$. Therefore $S_{\alpha, \sigma, \eps}r_o$ also converges weakly to $B$.
\qed
\bigskip

We now deduce Proposition~\ref{P:Donsker1} from Lemma~\ref{L:KDonsker}.
\medskip

\noindent
{\bf Proof of Proposition~\ref{P:Donsker1}.} As in (\ref{Xtau}), let $(r_o(T_i), T_i)_{i\geq 1}$ be the successive break points along $(r_o(n))_{n\geq 0}$. Note that
$\gamma_o(T_i)=r_o(T_i)$ for all $i\geq 1$. By the independence of $(\gamma_o, r_o)$ before and after a break point
conditional on the break point, we observe again that conditional on the first break point $(\gamma_o(T_1), T_1)=(r_o(T_1), T_1)$,
$$
(\gamma_o(T_{i-1}+j)-\gamma_o(T_{i-1}),\ r_o(T_{i-1}+j)-r_o(T_{i-1}))_{0\leq j\leq \tau_i}, \qquad i\geq 2,
$$
are i.i.d.\ and independent of $(\gamma_o(j), r_o(j))_{0\leq j\leq \tau_1}$. Suppose that $T_1<\infty$ a.s. Then by conditioning on $(r_o(T_1), T_1)$,
Lemma~\ref{L:KDonsker} implies that $S_{\alpha, \sigma, \eps}r_o$ converges weakly to a standard Brownian motion $B$
since $S_{\alpha, \sigma, \eps}(\sup_{t\in [0, T_1]} r_o(t)) \to 0$ and $S_{\alpha, \sigma, \eps} T_1 \to 0$ in probability as $\eps\downarrow 0$.
To conclude that $S_{\alpha, \sigma, \eps} \gamma_o$ and $S_{\alpha, \sigma, \eps} r_o$ converge to the same Brownian motion, it suffices to note that
\be\label{2tau2}
\max_{T_{i-1}\leq j\leq T_i} (r_o(j)-\gamma_o(j))  \leq 2\tau_i, \qquad i\geq 2,
\ee
are i.i.d.\ with all moments finite, since $r_o$ (resp.\ $\gamma_o$) can increase (resp.\ decrease) by at most $1$ each step, and $\tau_i$ has all moments finite by Lemma~\ref{L:KCLT}.

Now we show that $T_1<\infty$ a.s. We decompose the probability space by the value of $\gamma_o(0)=2x$ for $x\leq 0$. On the event $\gamma_o(0)=2x$
(i.e., $2x\in \Ki$ and $2i\not\in \Ki$ for $x+1\leq i\leq 0$), we note that $(r_o(T_1), T_1)$ is simply the first break point along the right boundary $r_{(2x, 0)}$ of the exploration cluster $C_{(2x,0)}$ after the finite open clusters at $(2i,0)$, $x+1\leq i\leq 0$, have all died out. Therefore on the event $\gamma_o(0)=2x$, $T_1<\infty$ a.s.\
by Lemma~\ref{L:KCLT}.
\qed
\bigskip

We conclude this section with an error bound on the approximation of $\gamma_o$ by $r_o$.
\bl{\bf [Approximation error]}\label{L:error} For each $L>0$ and $0<\delta<1$, there exists $C>0$ such that for all $\eps\in (0,1]$, we have
\be\label{error}
\P\Big(\sup_{t\in [0,\eps^{-1}L]}|r_o(t)- \gamma_o(t)| \geq \eps^{-\delta} \Big) \leq C \eps^{1/\delta},
\ee
and
\be\label{gap}
\P\Big(r_o(s)\neq \gamma_o(s) \ \forall\, s\in [t, t+\eps^{-\delta}] \ \mbox{for some } t\in [0,\eps^{-1}L] \Big) \leq C \eps^{1/\delta}.
\ee
\el
{\bf Proof.} As noted in the proof of Proposition~\ref{P:Donsker1}, conditional on the first break point $(r_o(T_1), T_1)$ along $r_o$,
$(r_o(j)-\gamma_o(j))_{0\leq j\leq T_1}$ and $r_o(T_{i-1}+j)-\gamma_o(T_{i-1}+j))_{0\leq j\leq \tau_i}$, $i\geq 2$, are independent, with
the latter forming an i.i.d.\ sequence.

Since $(\tau_i)_{i\geq 2}$ are i.i.d.\ and non-negative, $\frac{1}{n}\sum_{i=2}^{n+1} \tau_i$ satisfies a lower large deviation bound. In particular,
for any $c>\E[\tau_2]^{-1}$, there exist $C_1, C_2>0$ depending on $c$ and $L$ such that for all $\eps\in (0,1]$,
\be\label{ldp}
\begin{aligned}
\P(T_{\lceil c\eps^{-1}L\rceil+1}\leq \eps^{-1}L) & \leq
\P\Big(\frac{1}{\lceil c\eps^{-1}L\rceil} \sum_{i=2}^{\lceil c\eps^{-1}L\rceil+1}\tau_i \leq \frac{\eps^{-1}L}{\lceil c\eps^{-1}L\rceil} \Big)
\leq C_1 e^{-C_2\eps^{-1}L},
\end{aligned}
\ee
which decays faster than any power of $\eps$ as $\eps\downarrow 0$. Therefore to prove (\ref{error}), it suffices to show
\be\label{error2}
\P\Big(\sup_{t\in [0,T_{\lceil c\eps^{-1}L\rceil+1}]}|r_o(t)- \gamma_o(t)| \geq \eps^{-\delta} \Big) \leq C \eps^{1/\delta}.
\ee
We bound the supremum of $|r_o(t)- \gamma_o(t)|$ on $[0, T_1]$ and $[T_1, T_{\lceil c\eps^{-1}L\rceil +1}]$ separately.

Firstly,
$$
\begin{aligned}
\P\Big(\sup_{t\in [T_1,T_{\lceil c\eps^{-1}L\rceil+1}]}|r_o(t)- \gamma_o(t)| \geq \eps^{-\delta} \Big) &\leq  \sum_{i=2}^{\lceil c\eps^{-1}L\rceil+1} \P\big(\sup_{0\leq s\leq \tau_i} |r_o(T_{i-1}+s)-\gamma_o(T_{i-1}+s)|\geq \eps^{-\delta}\big) \\
&= \lceil c\eps^{-1}L\rceil \P\big(\sup_{0\leq s\leq \tau_2} |r_o(T_1+s)-\gamma_o(T_1+s)|\geq \eps^{-\delta}\big) \\
&\leq \lceil c\eps^{-1}L\rceil \eps^{k\delta}\, \E\Big[\sup_{0\leq s\leq \tau_2} |r_o(T_1+s)-\gamma_o(T_1+s)|^k\Big] \\
& \leq C \eps^{1/\delta}
\end{aligned}
$$
for some $C>0$ depending on $c$ and $L$, where we have applied the Markov inequality, chosen $k$ to be sufficiently large, and used the fact that $\sup_{0\leq s\leq \tau_2} |r_o(T_1+s)-\gamma_o(T_1+s)|$ has all moments finite as noted in (\ref{2tau2}).

Secondly, we note that by the same reasoning as for (\ref{2tau2}),
$$
\sup_{t\in [0, T_1]}|r_o(t)- \gamma_o(t)|\leq |\gamma_o(0)|+2T_1.
$$
Therefore to prove
$$
\P\big(\sup_{t\in [0, T_1]}|r_o(t)- \gamma_o(t)| \geq \eps^{-\delta} \big)\leq C\eps^{1/\delta}
$$
and thus deduce (\ref{error2}), it suffices to show that $|\gamma_o(0)|$ and $T_1$ have all moments finite.

Let $H_i$ be the time when the open cluster at $(-2i, 0)$ dies out. Then for any $k\geq 1$,
\be\label{gammaotail}
\P(|\gamma_o(0)| \geq 2k) = \P(H_i<\infty \mbox{ for all } 0\leq i < k) \leq C_3 e^{-C_4 k}
\ee
for some $C_3, C_4>0$ by results in \cite[Section 10]{D84}. Therefore $|\gamma_o(0)|$ has all moments finite.

On the other hand, given $\gamma_o(0)=-2x$ for some $x\geq 0$, $T_1$ is the first break point along $r_{(-2x, 0)}$ after the open clusters at
$(-2i, 0)$, $0\leq i \leq x-1$, have all died out. If we let $(r_{(-2x,0)}(\bar T_i), \bar T_i)$, $i\in\N$, denote the successive break points
along $r_{(-2x, 0)}$, with $\bar\tau_i=\bar T_i-\bar T_{i-1}$, then we have
\begin{eqnarray}
&& \!\!\!\! \P(T_1 \geq 2k)  = \sum_{x=0}^\infty \P(\gamma_o(0)=-2x, T_1\geq 2k) \nonumber\\
&\!\!\!\!\leq&\!\!\!\! \sum_{x=0}^\infty \P(\gamma_o(0)=-2x, T_1\geq 2k, k\leq \!\!\max_{0\leq i\leq x-1}\!\!\! H_i <\!\infty) + \!\!\sum_{x=0}^\infty \P(\gamma_o(0)=-2x, T_1\! \geq 2k, \!\max_{0\leq i\leq x-1}\!\!\! H_i <\! k) \nonumber\\
&\!\!\!\!\leq&\!\!\!\! \sum_{x=0}^\infty \sum_{i=0}^{x-1} \P(\gamma_0=-2x, k\leq H_i<\infty) + \sum_{x=0}^\infty \P(\gamma_o(0)=-2x, (-2x,0)\in\Ki, \max_{2\leq i\leq k} \bar \tau_i \geq k) \nonumber \\
&\!\!\!\!\leq&\!\!\!\! \sum_{x=0}^\infty x \P(\gamma_o(0)=-2x)^{\frac{1}{2}}\P(k\leq H_0<\! \infty)^{\frac{1}{2}} + \!\!\sum_{x=0}^\infty \P(\gamma_o(0)=-2x)^{\frac{1}{2}} \P((-2x,0)\in\Ki, \max_{2\leq i\leq k} \bar\tau_i \geq\! k)^{\frac{1}{2}} \nonumber \\
&\!\!\!\!\leq&\!\!\!\! C_5 \P(k\leq H_0< \infty)^{\frac{1}{2}} + C_6 \P(\max_{2\leq i\leq k} \tau_i \geq k | o\in \Ki)^{\frac{1}{2}}, \label{T1bound}
\end{eqnarray}
where we have applied Cauchy-Schwartz inequality and used (\ref{gammaotail}) to bound $\P(\gamma_o(0)=-2x)$. The first term in (\ref{T1bound}) decays exponentially in $k$, because
$$
\P(k\leq H_0< \infty)\leq C_7e^{-C_8k}
$$
for some $C_7, C_8>0$ by results in \cite[Section 12]{D84}. The second term in (\ref{T1bound}) decays faster than any power of $k$, because conditional on $o\in \Ki$, $(\tau_i)_{i\geq 2}$ are i.i.d.\ with all moments finite by Lemma~\ref{L:KCLT}. Therefore $T_1$ also has all moments finite. This concludes the proof of (\ref{error2}) and hence also (\ref{error}).

To prove (\ref{gap}), we note that $\gamma_o(T_i)=r_o(T_i)$ for all $i\in\N$. Therefore the event in (\ref{gap}) is contained in the event that there exists
some $i\in\N$ with $[T_{i-1}, T_i]\subset [0, \eps^{-1}L+\eps^{-\delta}]$ and $\tau_i \geq \eps^{-\delta}$. By the same large deviation bound as in (\ref{ldp}),
we may restrict to the event $\eps^{-1}L+\eps^{-\delta} \leq T_{\lceil c\eps^{-1}L\rceil+1}$. Then the probability in (\ref{gap}) can be bounded by
$$
\P(T_1 \geq \eps^{-\delta}) + \sum_{i=2}^{\lceil c\eps^{-1}L\rceil+1} \P(\tau_i\geq \eps^{-\delta}) = \P(T_1 \geq \eps^{-\delta}) + \lceil c\eps^{-1}L\rceil\P(\tau_2\geq \eps^{-\delta}) \leq C \eps^{1/\delta},
$$
again because $T_1$ and $\tau_2$ have all moments finite.
\qed

\section{Convergence to Coalescing Brownian Motions}\label{S:I}
In this section, we show that different exploration clusters evolve independently before they intersect, and when two exploration clusters intersect, they coalesce in most cases. As a consequence, we prove that after applying the shearing and diffusive scaling map $S_{\alpha, \sigma, \eps}$, a finite number of exploration clusters $C_{z_i}$, $1\leq i\leq k$, converge in distribution to coalescing Brownian motions, which implies the convergence criterion (I) for $(S_{\alpha, \sigma,\eps}\overline{\Gamma})_{\eps\in (0,1)}$.

\subsection{Convergence of a pair of exploration clusters}
\begin{figure}[tp] 
\begin{center}
\includegraphics[height=19cm]{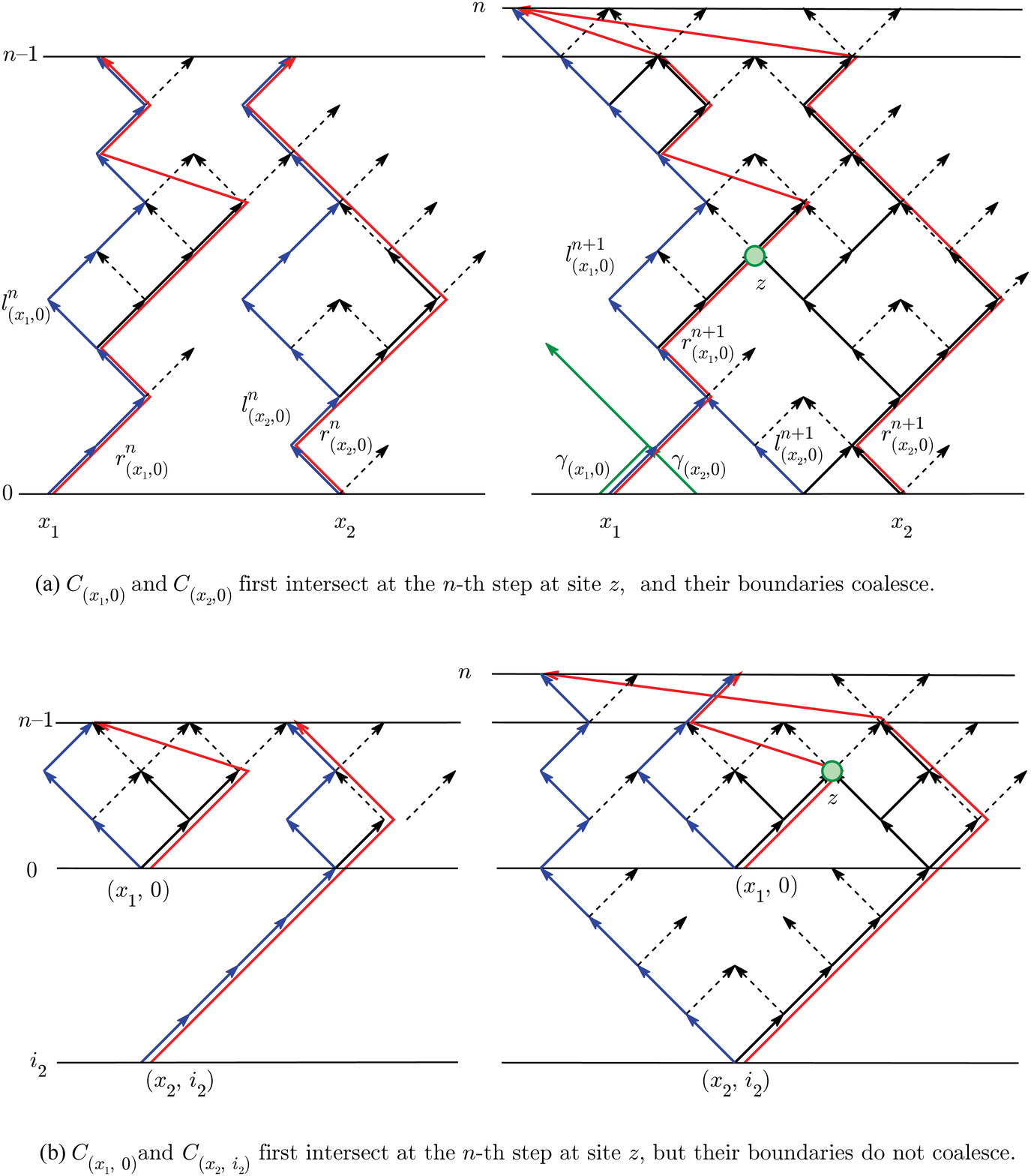}
\caption{When two exploration clusters first intersect.}
\label{fig:coalescence}
\end{center}
\end{figure}

By construction, two exploration clusters evolve independently until the first time they intersect, i.e., when they share a common explored edge.
We first show that two exploration clusters starting at the same time must coalesce when they intersect. Complications arise when the two clusters start at different times. See Figure~\ref{fig:coalescence}.

\bl{\bf [Coalescence of exploration clusters]}\label{L:coalesce}
Let $z_1=(x_1, 0), z_2=(x_2, 0)\in\Z^2_{\rm even}$ with $x_1< x_2$. Let $C_{z_i}$ and $(l^n_{z_i}, r^n_{z_i})_{n\geq 0}$, $i=1,2$, be the
respective exploration cluster processes and their associated boundary processes. Let $(r_{z_i}(n))_{n\geq 0}$, $i=1,2$, be the infinite right boundaries defined as in Lemma~\ref{L:lrprop}. Define
$$
\begin{aligned}
\kappa_{rl} & :=\min\{n\geq 0: l^n_{z_2}(i)\leq r_{z_1}(i)   \mbox{ for some } 0\leq i\leq n\}, \\
\kappa_{rr} & := \min\{n\geq 0: r_{z_2}(n) \leq r_{z_1}(n)\}, \\
\kappa_{\gamma\gamma} & :=\min\{n\geq 0: \gamma_{z_2}(n)\leq \gamma_{z_1}(n)\}.
\end{aligned}
$$
Then
\begin{itemize}
\item[\rm (i)] $\kappa_{rr}=\kappa_{rl}$ and $r_{z_1}(n)=r_{z_2}(n)$ for all $n\geq \kappa_{rr}$.

\item[\rm (ii)] $l^n_{z_1}(\cdot)=l^n_{z_2}(\cdot)$ on $[\kappa_{rr}, n]$ for all $n\geq \kappa_{rr}$.

\item[\rm (iii)] $\kappa_{\gamma\gamma} =\min\{n\geq 0: \gamma_{z_2}(n)\leq r_{z_1}(n)\}$, $\kappa_{\gamma\gamma}\leq \kappa_{rr}$, and $\gamma_{z_1}(n)=\gamma_{z_2}(n)$ for all $n\geq \kappa_{\gamma\gamma}$.
\end{itemize}
If $z_2=(x_2, i_2)\in\Z^2_{\rm even}$ with $i_2\neq 0$, then the above statements hold on the event $r_{z_1}(0\vee i_2) \leq \gamma_{z_2}(0\vee i_2)$, with $0$ replaced by $0\vee i_2$ in the definition of $\kappa_{rl}$, $\kappa_{rr}$ and $\kappa_{\gamma\gamma}$. By symmetry, analogous statements hold on the event $r_{z_2}(0\vee i_2)\leq \gamma_{z_1}(0\vee i_2)$.
\el

\br{\rm Lemma~\ref{L:coalesce} shows that if $z_1=(x_1, i_1)$ and $z_2=(x_2,i_2)$ satisfy $i_1=i_2$, or $r_{z_1}(i_1\vee i_2) \leq \gamma_{z_2}(i_1\vee i_2)$, or $r_{z_2}(i_1\vee i_2) \leq \gamma_{z_1}(i_1\vee i_2)$, then the time when $C_{z_1}(\cdot)$ and $C_{z_2}(\cdot)$ first intersect is also the time when $r_{z_1}$
and $r_{z_2}$ coalesce, and $r_{z_1}$ and $r_{z_2}$ cannot intersect or cross each other without coalescing. This picture may fail when the above conditions
are violated: one such scenario is illustrated in Figure~\ref{fig:coalescence}~(b) where $l^n_{z_2}(0)<l^n_{z_1}(0)= r_{z_1}(0)<r_{z_2}(0)$. }
\er

\bigskip
\noindent
{\bf Proof of Lemma~\ref{L:coalesce}.} For (i)--(ii): Since $l^n_{z_2}(\cdot) \leq r_{z_2}(\cdot)$ on $[0, n]$ for all $n\in\N$, clearly $r_{z_1}(\cdot)<r_{z_2}(\cdot)$
on $[0, \kappa_{rl}-1]$ and hence $\kappa_{rr}\geq \kappa_{rl}$. Now let $n\geq \kappa_{rl}$.

If $l_{z_2}^n(0)< x_1$, then because $l_{z_i}^n$, $i=1,2$, is the rightmost open path connecting $(-\infty, x_i]\times\{0\}$ to $\Z\times\{n\}$, we must have $l_{z_1}^n=l_{z_2}^n$, which by Lemma~\ref{L:lrprop}~(iii) implies that $r_{z_1}(n)=r_{z_2}(n)$.

Now suppose that $x_1\leq l_{z_2}^n(0)$. By the definition of $\kappa_{rl}$ and Lemma~\ref{L:lrprop}~(iv) that $l^n_{z_2}(\cdot)$ decreases in $n$, we must have $l_{z_2}^n(i)\leq r_{z_1}(i)$ for some $0\leq i\leq \kappa_{rl}$. Since $r_{z_1}$ can only increase by at most one at each step while $l_{z_2}^n$ is a nearest-neighbor path, $l_{z_2}^n$ and $r_{z_1}$ must coincide at some time in $\{0,1,\cdots, \kappa_{rl}\}$. Let $j$ be the first such time. By the definition of $r_{z_1}$, the rightmost open path $\pi$ connecting $(-\infty, x_1]\times\{0\}$ to $\Z\times\{j\}$ satisfies $\pi(j)=r_{z_1}(j)=l_{z_2}^n(j)$. By concatenating first $\pi$ and then $l_{z_2}^n$ at time $j$, we obtain an open path connecting $(-\infty, x_1]\times\{0\}$ to $\Z\times\{n\}$, which coincides with $l^n_{z_2}(\cdot)$ on $[j, n]\supset [\kappa_{rl}, n]$, lies on the right of $l^n_{z_1}(\cdot)$ on $[0, j]$ by our choice of $\pi$, and lies on the left of $l^n_{z_1}(\cdot)$ on $[0, n]$ by the definition of $l^n_{z_1}$. Therefore $l^n_{z_2}(\cdot)\leq l^n_{z_1}(\cdot)$ on $[j,n]$ and $l^n_{z_1}(j)=l^n_{z_2}(j)$. We must also have $l^n_{z_1}(\cdot)\leq l^n_{z_2}(\cdot)$ on $[j,n]$, because concatenating first $l^n_{z_2}$ and then $l^n_{z_1}$ at time $j$ produces an open path connecting $(-\infty, x_2]\times\{0\}$ to $\Z\times \{n\}$ and $l^n_{z_2}$ is the rightmost such open path. Therefore $l^n_{z_1}(\cdot)=l^n_{z_2}(\cdot)$ on $[j, n]\supset [\kappa_{rl}, n]$, which further implies $r_{z_1}(n)=r_{z_2}(n)$ by Lemma~\ref{L:lrprop}~(iii). Combining the conclusions above then gives (i) and (ii).

For (iii), the argument is the same. We distinguish between the cases $\gamma_{z_2}(0)< x_1$ and $x_1\leq \gamma_{z_2}(0)$. In the first case we deduce
$\gamma_{z_1}=\gamma_{z_2}$ from their definition. In the second case, we find the first time $j$ when $r_{z_1}$ and $\gamma_{z_2}$ coincide. By concatenating paths
as done above and using the definitions of $\gamma_{z_1}$ and $\gamma_{z_2}$, we deduce that $\gamma_{z_1}(\cdot)$ and $\gamma_{z_2}(\cdot)$ must coincide on $[j, \infty)$. Since $\gamma_{z_1}\leq r_{z_1}$, $\gamma_{z_1}$ and $\gamma_{z_2}$ cannot intersect earlier. Therefore $j=\kappa_{\gamma\gamma}$. Since $\gamma_{z_2}\leq r_{z_2}$, we have $\kappa_{\gamma\gamma}\leq\kappa_{rr}$.

In the proof above, the fact that $z_1$ and $z_2$ are at the same time is used when we tried to prove $l^n_{z_1}=l^n_{z_2}$ in the scenario $l^n_{z_2}(0)< x_1=r_{z_1}(0)$ for some $n\geq \kappa_{rl}$. When $z_1$ and $z_2$ are at different times, $l^n_{z_1}$ and $l^n_{z_2}$ are in general not comparable and the same reasoning does not apply.
However, the condition $r_{z_1}(0\vee i_2)\leq \gamma_{z_2}(0\vee i_2)$ rules out such a scenario and guarantees that $r_{z_1}(0\vee i_2)\leq l^n_{z_2}(0\vee i_2)$ for all $n\geq 0\vee i_2$. The concatenation arguments we had for $i_2=0$ then applies without change.
\qed
\bigskip

We now formulate the convergence of a pair of exploration clusters to coalescing Brownian motions, which extends Proposition~\ref{P:Donsker1} for a single exploration cluster.

\bprop\label{P:Donsker2}{\bf [Convergence of a pair of exploration clusters]} Let $\alpha$ and $\sigma$ be as in (\ref{alphasigma}). For $i=1, 2$, let
$z^\eps_i = (x^\eps_i, n^\eps_i)\in\Z^2_{\rm even}$ be such that $S_{\alpha, \sigma, \eps}z^\eps_i \to z_i=(x_i, t_i)\in\R^2$ as $\eps\downarrow 0$. Let
\be\label{kapparg}
\begin{aligned}
\kappa^\eps_{rr} & := \min\{n\in\Z: r_{z^\eps_1}(m) = r_{z^\eps_2}(m)\ \forall\, m\geq n\} \\
\kappa^\eps_{\gamma\gamma} & := \min\{n\in\Z: \gamma_{z^\eps_1}(n) = \gamma_{z^\eps_2}(n) \ \forall\, m\geq n\}.
\end{aligned}
\ee
Let $B_1$ and $B_2$ be two coalescing Brownian motions starting at respectively $z_1$ and $z_2$, and let
\be\label{kappa}
\kappa:=\inf\{t\in\R: B_1(t)=B_2(t)\}.
\ee
Then as $\eps\downarrow 0$,
\be\label{conv2}
S_{\alpha, \sigma, \eps}(\gamma_{z^\eps_1}, r_{z^\eps_1}, \gamma_{z^\eps_2}, r_{z^\eps_2}, \kappa^\eps_{\gamma\gamma}, \kappa^\eps_{rr})
\stackrel{\rm dist}{\Longrightarrow} (B_1, B_1, B_2, B_2, \kappa, \kappa)
\ee
as random variables taking values in the product space $\Pi^4\times [-\infty,\infty]^2$.
\eprop
{\bf Proof.} Our proof strategy is similar to the proof that two coalescing random walks converge in distribution to two coalescing Brownian motions. We first recall the argument in that context to serve as a guide. Start with two independent random walks, which converge in distribution to two independent Brownian motions by Donsker's invariance principle. Using Skorohod's representation theorem for weak convergence~\cite[Theorem 6.7]{B99}, we can use coupling to turn such a convergence into almost sure convergence in path space. Next we observe that coalescing random walk paths can be constructed almost surely from two independent random walk paths by forcing the second walk to follow the first walk from the moment they meet, and the same deterministic operation applied to two independent Brownian motions gives a construction of two coalescing Brownian motions. It is not difficult to show that under the coupling given by Skorohod's representation theorem, almost surely the time the two independent walks meet converges to the time the two independent Brownian motions meet, which then implies that the two coalescing random walks constructed above also converge almost surely to the two coalescing Brownian motions.

To apply the above argument to our context, we first construct the two exploration clusters $(C_{z^\eps_1}, C_{z^\eps_2})$ from two independent percolation exploration clusters $C_{z^\eps_1}^{[1]}$ and $C_{z^\eps_2}^{[2]}$, equally distributed with $C_{z^\eps_1}$ and $C_{z^\eps_2}$ respectively. More precisely, let
\be\label{Omega}
\Omega^{\eps}_{[i]}:=\{\omega^{\eps,\pm}_{[i],u} : u\in\Z^2_{\rm even}\}, \quad i=1,2,
\ee
be two independent percolation configurations, where $\omega^{\eps,\pm}_{[i],u}$ are i.i.d.\ Bernoulli random variables with parameter $p$, and $\omega^{\eps,+}_{[i], u}$, resp.\ $\omega^{\eps,-}_{[i], u}$, equals $1$ if
the up-right, resp.\ up-left, edge from $u$ is open in the $i$-th percolation configuration and equals $0$ otherwise. Let $C^{[i]}_{z^\eps_i}$ be the exploration cluster starting at $z^\eps_i$ constructed from the percolation configuration $\Omega^\eps_{[i]}$, and let $\gamma_{z^\eps_i}^{[i]}$, $(l_{z^\eps_i}^{[i],n})_{n\geq n^\eps_i}$ and
$r_{z^\eps_i}^{[i]}$ be the associated rightmost infinite open path and left and right exploration cluster boundaries. Clearly $C_{z^\eps_1}^{[1]}$ and $C_{z^\eps_2}^{[2]}$ are independent, and equally distributed with $C_{z^\eps_1}$ and $C_{z^\eps_2}$ respectively.

Next we show how to construct $C_{z^\eps_1}$ and $C_{z^\eps_2}$ almost surely from $C_{z^\eps_1}^{[1]}$ and $C_{z^\eps_2}^{[2]}$ (more precisely, from the underlying percolation configurations $\Omega^{\eps}_{[1]}$ and $\Omega^{\eps}_{[2]}$), akin to the almost sure construction of two coalescing random walk paths from two independent random walk paths. First we set $C_{z^\eps_1}:=C^{[1]}_{z^\eps_1}$ so that $C_{z^\eps_1}$ is constructed by exploring the status of edges in $\Omega^\eps_{[1]}$. Next we construct $C_{z^\eps_2}$ by exploring the status of edges in $\Omega^\eps_{[2]}$ one by one until we first encounter an edge whose status in $\Omega^{\eps}_{[1]}$ has already been explored in the construction of $C_{z^\eps_1}$. From this step onward, we continue the exploration construction of $C_{z^\eps_2}$ by using the status of edges which have either already been explored so far in the construction of $C_{z^\eps_1}$ and $C_{z^\eps_2}$, or if a previously unexplored edge is encountered, we just look up its status in $\Omega^\eps_{[1]}$. Because the status of edges are discovered in a Markovian way, $(C_{z^\eps_1}, C_{z^\eps_2})$ constructed this way has the right distribution. Note that $C_{z^\eps_1}=C^{[1]}_{z^\eps_1}$, and $C_{z^\eps_2}(\cdot) = C^{[2]}_{z^\eps_2}(\cdot)$ on $[n^\eps_2, \iota^{[12],\eps}-1]$, where $\iota^{[12],\eps}$ is the first time $n$ when $C^{[2]}_{z^\eps_2}(n)$ encounters an edge which has already been explored in the construction of $C_{z^\eps_1}=C^{[1]}_{z^\eps_1}$. In particular,
\be\label{r12iota}
r_{z^\eps_1} = r^{[1]}_{z^\eps_1}, \quad \mbox{and}\quad r_{z^\eps_2}(\cdot) = r^{[2]}_{z^\eps_2}(\cdot) \mbox{ on } [n^\eps_2, \iota^{[12],\eps}-1].
\ee
If $r_{z^\eps_1}$ and $r_{z^\eps_2}$ coalesce at time $\iota^{[12],\eps}$, then the analogy with the proof of convergence of coalescing random walks to coalescing Brownian motions will be complete. However this is not always true, and one such case has been explained in the remark after Lemma~\ref{L:coalesce} and illustrated in Figure~\ref{fig:coalescence}~(b). Fortunately, such events are rare and can be controlled using Lemma~\ref{L:error}, which we show next.

Note that $\{S_{\alpha, \sigma, \eps}(\gamma_{z^\eps_1}, r_{z^\eps_1}, \gamma_{z^\eps_2}, r_{z^\eps_2}, \kappa^\eps_{\gamma\gamma}, \kappa^\eps_{rr})\}_{\eps>0}$ is a tight family of $\Pi^4\times [-\infty,\infty]^2$-valued random variables because
$S_{\alpha, \sigma, \eps}(\gamma_{z^\eps_i}, r_{z^\eps_i}) \stackrel{\rm dist}{\Longrightarrow} (B_i, B_i)$, $i=1,2$, by Proposition~\ref{P:Donsker1}. Therefore it suffices to verify (\ref{conv2}) along any weakly convergent subsequence. By Lemma~\ref{L:error}, for each $L\in\N$ and $0<\delta<1$, there exists $C_{\delta, L}$ such that
$$
\max_{i=1, 2}\ \P\Big(\sup_{t\in [\eps n^\eps_i, L]} |S_{\alpha, \sigma, \eps}r_{z^\eps_i}(t) - S_{\alpha, \sigma, \eps}\gamma_{z^\eps_i}(t)| \geq \eps^{\frac{1}{2}-\delta}\Big) \leq C_{\delta, L} \eps^{1/\delta},
$$
$$
\max_{i=1,2}\ \P\big(r_{z^\eps_i}(s/\eps)\neq \gamma_{z^\eps_i}(s/\eps)
\ \forall\, s\in [t, t+\eps^{1-\delta}] \ \mbox{for some } t\in [\eps n^\eps_i, L] \big) \leq C_{\delta, L} \eps^{1/\delta}.
$$
Let $\delta=\frac{1}{4}$. By going to a further subsequence if necessary, it suffices to verify (\ref{conv2}) along any weakly convergent subsequence indexed by $(\eps_m)_{m\in\N}$ with $\eps_m\downarrow 0$, such that for all $L\in\N$,
\be\label{errorBC}
\sum_{m=1}^\infty \max_{i=1, 2}\ \P\Big(\sup_{t\in [\eps_m n^{\eps_m}_i, L]} |S_{\alpha, \sigma, \eps_m}r_{z^{\eps_m}_i}(t) - S_{\alpha, \sigma, \eps_m}\gamma_{z^{\eps_m}_i}(t)| \geq \eps_m^{1/4}\Big) <\infty,
\ee
\be\label{gapBC}
\sum_{m=1}^\infty  \max_{i=1,2}\ \P\big(r_{z^{\eps_m}_i}(s/\eps_m)\neq \gamma_{z^{\eps_m}_i}(s/\eps_m)
\ \forall\, s\in [t, t+\eps_m^{3/4}] \ \mbox{for some } t\in [\eps_m n^{\eps_m}_i, L] \big)<\infty.
\ee
By Borel-Cantelli, almost surely, the events in (\ref{errorBC}) and (\ref{gapBC}) happen only a finite number of times regardless of how the percolation configurations $(\Omega^{\eps_m}_{[1]}, \Omega^{\eps_m}_{[2]})_{m\in\N}$ are coupled. In words, (\ref{errorBC}) implies that $S_{\alpha, \sigma, \eps_m}r_{z^{\eps_m}_i}$ and $S_{\alpha, \sigma, \eps_m}\gamma_{z^{\eps_m}_i}$ are almost surely close as $m\to\infty$, while (\ref{gapBC}) implies that the maximum gap between successive regeneration times along $S_{\alpha, \sigma, \eps_m}r_{z^{\eps_m}_i}$ and $S_{\alpha, \sigma, \eps_m}\gamma_{z^{\eps_m}_i}$ almost surely tends to $0$ as $m\to\infty$. We will need both properties later. From now on we work with such a sequence of $(\eps_m)_{m\in\N}$.

By Proposition~\ref{P:Donsker1}, as $m\to\infty$,
\be\label{preskorohod}
S_{\alpha, \sigma, \eps_m}(\gamma^{[1]}_{z^{\eps_m}_1}, r^{[1]}_{z^{\eps_m}_1}, \gamma^{[2]}_{z^{\eps_m}_2}, r^{[2]}_{z^{\eps_m}_2}) \longrightarrow (B^{[1]}_1, B^{[1]}_1, B^{[2]}_2, B^{[2]}_2)
\ee
in distribution for two independent Brownian motions $B^{[1]}_1$ and $B^{[2]}_2$, starting respectively at $z_1$ and $z_2$. By Skorohod's representation theorem, we can couple the sequence of random variables $\{S_{\alpha, \sigma, \eps_m}(\gamma^{[1]}_{z^{\eps_m}_1}, r^{[1]}_{z^{\eps_m}_1}, \gamma^{[2]}_{z^{\eps_m}_2}, r^{[2]}_{z^{\eps_m}_2})\}_{m\in\N}$ and $(B^{[1]}_1$, $B^{[2]}_2)$ on the same probability space such that the convergence in (\ref{preskorohod}) becomes almost sure. Furthermore such a coupling can be extended to a coupling of the underlying sequence of percolation configurations
$(\Omega^{\eps_m}_{[1]}, \Omega^{\eps_m}_{[2]})_{m\in\N}$ by sampling $(\Omega^{\eps_m}_{[1]}, \Omega^{\eps_m}_{[2]})$, $m\in\N$, independently conditional on the realization of $(\gamma^{[1]}_{z^{\eps_m}_1}, r^{[1]}_{z^{\eps_m}_1}, \gamma^{[2]}_{z^{\eps_m}_2}, r^{[2]}_{z^{\eps_m}_2})$, $m\in\N$. Let us assume such a coupling from now on. We will show that the convergence in (\ref{conv2}) in fact takes place almost surely, similar in spirit to the proof that two coalescing random walks converge to two coalescing Brownian motions almost surely, once the independent random walks and Brownian motions used to construct the coalescing systems are coupled properly.

Let
$$
\kappa^{[12]}:=\inf\{t\in\R: B^{[1]}_1(t) = B^{[2]}_2(t)\}.
$$
Define $B_1:=B^{[1]}_1$, and $B_2(\cdot):=B^{[2]}_2(\cdot)$ on $[t_2, \kappa^{[12]}]$ and
$B_2(\cdot) = B^{[1]}_1(\cdot)$ on $[\kappa^{[12]}, \infty]$. Then $(B_1, B_2)$ is distributed as a pair of coalescing Brownian motions starting respectively
at $z_1$ and $z_2$, and $\kappa^{[12]}=\kappa$ as defined in (\ref{kappa}). Let $C_{z^{\eps_m}_i}$, $\gamma_{z^{\eps_m}_i}$, $(l_{z^{\eps_m}_i}^{[i],n})_{n\geq n^{\eps_m}_i}$ and $r_{z^{\eps_m}_i}$, $i=1,2$, be constructed from the percolation configurations $(\Omega^{\eps_m}_{[1]}, \Omega^{\eps_m}_{[2]})$ as before. By construction, $(\gamma_{z^{\eps_m}_1}, r_{z^{\eps_m}_1})=(\gamma^{[1]}_{z^{\eps_m}_1}, r^{[1]}_{z^{\eps_m}_1})$. Therefore by (\ref{preskorohod}) and our coupling,
$$
S_{\alpha, \sigma, \eps_m}(\gamma^{[1]}_{z^{\eps_m}_1}, r^{[1]}_{z^{\eps_m}_1}) \to (B_1, B_1) \qquad {\rm a.s.}
$$
Assume first $z_1\neq z_2$. Then a.s.\ either (1) $y_1:=B_1(t_1\vee t_2) < y_2:=B_2(t_1\vee t_2)$, or (2) $y_2<y_1$.

In case (1), define
$$
\begin{aligned}
\kappa^{[12], \eps_m}_{r\gamma} & := \min\{n\geq n^{\eps_m}_1 \vee n^{\eps_m}_2: \gamma^{[2]}_{z^{\eps_m}_2}(n) \leq r^{[1]}_{z^{\eps_m}_1}(n)\}, \\
\kappa^{[12], \eps_m}_{rr} & := \min\{n\geq n^{\eps_m}_1 \vee n^{\eps_m}_2: r^{[2]}_{z^{\eps_m}_2}(n) \leq r^{[1]}_{z^{\eps_m}_1}(n)\}.
\end{aligned}
$$
Because the left boundary of the exploration cluster $C^{[2]}_{z^{\eps_m}_2}(n)$ is bounded between $\gamma^{[2]}_{z^{\eps_m}_2}$ and $r^{[2]}_{z^{\eps_m}_2}$, for $m$ sufficiently large, the time $\iota^{[12],\eps_m}$ when $C^{[2]}_{z^{\eps_m}_2}(n)$ first intersects $C_{z^\eps_1}=C^{[1]}_{z^{\eps_m}_1}$ satisfies
$$
\kappa^{[12], \eps_m}_{r\gamma} \leq \iota^{[12],\eps_m} \leq \kappa^{[12], \eps_m}_{rr}.
$$
Since the a.s.\ convergence in (\ref{preskorohod}) induced by our coupling implies that $(\eps_m \kappa^{[12], \eps_m}_{r\gamma}, \eps_m \kappa^{[12], \eps_m}_{rr})\to (\kappa, \kappa)$, we must have $\eps_m\iota^{[12],\eps_m}\to \kappa$ as well. By (\ref{r12iota}), $r_{z^{\eps_m}_2}(\cdot)=r^{[2]}_{z^{\eps_m}_2}(\cdot)$ on $[n^{\eps_m}_2, \iota^{[12], \eps_m}-1]$, and hence by (\ref{preskorohod}), $S_{\alpha, \sigma, \eps_m}r_{z^{\eps_m}_2}\to B_2$ uniformly on $[t_2, \kappa]$. To draw the same conclusion for $S_{\alpha, \sigma, \eps_m}\gamma_{z^{\eps_m}_2}$, we can no longer appeal to (\ref{preskorohod}) because there is no analogue of (\ref{r12iota}) that reduces $\gamma_{z^{\eps_m}_2}$ to $\gamma^{[2]}_{z^{\eps_m}_2}$. Instead we can use (\ref{errorBC}) and Borel-Cantelli to conclude that $S_{\alpha, \sigma, \eps_m}\gamma_{z^{\eps_m}_2}\to B_2$ uniformly on $[t_2, \kappa]$. It then follows that for all $m$ sufficiently large,
$$
r_{z^{\eps_m}_1}(n^{\eps_m}_1 \vee n^{\eps_m}_2) <\gamma_{z^{\eps_m}_2}(n^{\eps_m}_1 \vee n^{\eps_m}_2),
$$
which allows us to apply Lemma~\ref{L:coalesce} (note that $\iota^{[12],\eps_m}$ here is exactly $\kappa_{rl}$ in Lemma~\ref{L:coalesce}) to conclude that $r_{z^{\eps_m}_2}(\cdot) = r_{z^{\eps_m}_1}(\cdot)$ on $[\iota^{[12],\eps_m}, \infty)$
and $\kappa^{\eps_m}_{rr}=\iota^{[12], \eps_m}$. Therefore $S_{\alpha, \sigma, \eps_m} r_{z^{\eps_m}_2} \to B_2$ and $\eps\kappa^{\eps_m}_{rr}\to \kappa$. Again by
(\ref{errorBC}) and Borel-Cantelli, we also have $S_{\alpha, \sigma, \eps_m} \gamma_{z^{\eps_m}_2} \to B_2$. Since $\kappa^{\eps_m}_{\gamma\gamma} \leq \kappa^{\eps_m}_{rr}$ by Lemma~\ref{L:coalesce} and $S_{\alpha, \sigma, \eps_m}(\gamma_{z^{\eps_m}_1}, \gamma_{z^{\eps_m}_2})\to (B_1, B_2)$, we must also have $\eps\kappa^{\eps_m}_{\gamma\gamma}\to\kappa$. This proves (\ref{conv2}) with a.s.\ convergence under our coupling in case (1).

In case (2), define
$$
\kappa^{[12], \eps_m}_{r\gamma} := \min\{n\geq n^{\eps_m}_1 \vee n^{\eps_m}_2: \gamma^{[1]}_{z^{\eps_m}_1}(n) \leq r^{[2]}_{z^{\eps_m}_2}(n)\}.
$$
By (\ref{preskorohod}), $\eps_m\kappa^{[12], \eps_m}_{r\gamma}\to \kappa$.
Note that for all $m$ sufficiently large, we have
$$
\gamma^{[1]}_{z^{\eps_m}_1}(\kappa^{[12],\eps_m}_{r\gamma}) = r^{[2]}_{z^{\eps_m}_2}(\kappa^{[12],\eps_m}_{r\gamma}) \qquad \mbox{and}\qquad \iota^{[12], \eps_m}-1= \kappa^{[12], \eps_m}_{r\gamma}.
$$
These facts imply that $r_{z^{\eps_m}_2}(\cdot) = r^{[2]}_{z^{\eps_m}_2}(\cdot)$ on $[n^{\eps_m}_2, \kappa^{[12],\eps_m}_{r\gamma}]$. Furthermore, $(r_{z^{\eps_m}_2}(\kappa^{[12],\eps_m}_{r\gamma}), \kappa^{[12],\eps_m}_{r\gamma})$ is a break point along
$r_{z^{\eps_m}_2}$, and hence $\kappa^{[12],\eps_m}_{r\gamma}= \kappa^{\eps_m}_{\gamma\gamma}$ and $\eps_m \kappa^{\eps_m}_{\gamma\gamma}\to \kappa$. On $[\kappa^{\eps_m}_{\gamma\gamma},\infty)$, $r_{z^{\eps_m}_2}(\cdot)$ is bounded between $\gamma_{z^{\eps_m}_1}$ and $r_{z^{\eps_m}_1}$. Therefore by (\ref{preskorohod}), $S_{\alpha, \sigma, \eps_m}r_{z^{\eps_m}_2} \to B_2$. Applying (\ref{errorBC}) again with Borel-Cantelli gives $S_{\alpha, \sigma, \eps_m}\gamma_{z^{\eps_m}_2} \to B_2$. Finally we note that $\kappa^{\eps_m}_{rr}$ is bounded between $\kappa^{\eps_m}_{\gamma\gamma}$ and the first time after $\kappa^{\eps_m}_{\gamma\gamma}$ when $r_{z^{\eps_m}_1}(\cdot)=\gamma_{z^{\eps_m}_1}(\cdot)$. Then by (\ref{gapBC}) and Borel-Cantelli, we must have $\eps_m \kappa^{\eps_m}_{rr}\to \kappa$ as well. This proves (\ref{conv2}) with a.s.\ convergence under our coupling in case (2), and completes the proof of (\ref{conv2}) for $z_1\neq z_2$.

Lastly we treat the case $z_1=z_2$, which we may assume to be $o$ without loss of generality. For each $m\geq 3$, let $z^\eps_m\in\Z^2_{\rm even}$ be such that
$S_{\alpha, \sigma, \eps}z^\eps_m \to z_m=(1/m, 0)$ as $\eps\downarrow 0$. For $1\leq i< j$, let $\kappa^{\eps, ij}_{rr}$ and $\kappa^{\eps, ij}_{\gamma\gamma}$ be defined for $C_{z^{\eps}_i}$ and $C_{z^\eps_j}$ in the same way as in (\ref{kapparg}). Then
$$
\{S_{\alpha, \sigma, \eps} \big((r_{z^\eps_i}, \gamma_{z^\eps_i})_{i\in\N}, (\kappa^{\eps,ij}_{rr}, \kappa^{\eps, ij}_{\gamma\gamma})_{1\leq i< j}\big)\}_{\eps \in (0,1)}
$$
is a tight family of random variables taking values in $\Pi^\N\times [-\infty, \infty]^\N$ under the product topology. By (\ref{conv2}) proved earlier for $z_1\neq z_2$,
any weak limit must be of the form
$$
((W_i, W_i)_{i\in\N}, (\kappa^{12}_{rr}, \kappa^{12}_{\gamma\gamma}), (\kappa^{ij}, \kappa^{ij})_{1\leq i<j,\, j\geq 3})
$$
where for each $(i,j)$ with $1\leq i<j$ and $j\geq 3$, $(W_i, W_j, \kappa^{ij})$ is distributed as a pair of coalescing Brownian motions starting respectively at
$z_i$ and $z_j$, and $\kappa^{ij} =\inf\{t\geq 0: W_i(t)=W_j(t)\}$. Note that
$$
\kappa^{\eps, 12}_{rr} \leq \max\{ \kappa^{\eps, 1m}_{rr}, \kappa^{\eps, 2m}_{rr}\} \quad\mbox{and}\quad \kappa^{\eps, 12}_{\gamma\gamma} \leq \max\{ \kappa^{\eps, 1m}_{\gamma\gamma}, \kappa^{\eps, 2m}_{\gamma\gamma}\}
$$
for all $m\geq 3$, and $\kappa^{1m}\stackrel{\rm dist}{=}\kappa^{2m}$ converges in distribution to $0$ as $m\to \infty$. Therefore we must have $\kappa^{12}_{rr}=\kappa^{12}_{\gamma\gamma}=0$ a.s. It then follows that $W_1=W_2$ a.s. This concludes the proof of (\ref{conv2}) for $z_1=z_2$.
\qed

\br{\rm
Prop.~\ref{P:Donsker1}, Lemma~\ref{L:coalesce}, and the construction of two exploration clusters from two independent copies in the proof of
Prop.~\ref{P:Donsker2} show that two exploration clusters starting at the same time must coalesce a.s.\ in finite time. The same is true if two exploration
clusters start at different times, since each vertex can only reach a finite number of vertices by open path at any later time. This recovers the main result in~\cite{WZ08}, that any two paths in $\Gamma$ must coalesce a.s.\ in finite time.
}
\er

\subsection{Convergence of multiple exploration clusters}
We now extend Proposition~\ref{P:Donsker2} by establishing the convergence of a finite number of exploration clusters to coalescing Brownian motions, which implies that the convergence criterion (I) in Theorem~\ref{T:convct1} holds for $S_{\alpha, \sigma, \eps}\overline{\Gamma}$ as $\eps\downarrow 0$.

\bprop\label{P:Donsker3}{\bf [Convergence of multiple exploration clusters]} Let $k\in\N$. For $1\leq i\leq k$, let
$z^\eps_i = (x^\eps_i, n^\eps_i)\in\Z^2_{\rm even}$ be such that $S_{\alpha, \sigma, \eps}z^\eps_i \to z_i=(x_i, t_i)\in\R^2$ as $\eps\downarrow 0$.
For $1\leq i<j\leq k$, let $\kappa^{\eps, ij}_{rr}$ and $\kappa^{\eps, ij}_{\gamma\gamma}$ be defined for the exploration clusters $C_{z^\eps_i}$ and $C_{z^\eps_j}$
as in (\ref{kapparg}). Let $(B_1, \cdots, B_k)$ be coalescing Brownian motions starting respectively at $(z_1,\cdots, z_k)$, and let
$\kappa^{ij}$ be the time when $B_j$ and $B_j$ coalesce. Then as $\eps\downarrow 0$,
\be\label{conv3}
S_{\alpha, \sigma, \eps}\big((\gamma_{z^\eps_i}, r_{z^\eps_i})_{1\leq i\leq k}, (\kappa_{\gamma\gamma}^{\eps, ij}, \kappa_{rr}^{\eps, ij})_{1\leq i<j\leq k}\big)
\stackrel{\rm dist}{\Longrightarrow} \big((B_i, B_i)_{1\leq i\leq k}, (\kappa^{ij}, \kappa^{ij})_{1\leq i<j\leq k}\big)
\ee
as random variables taking values in the product space $\Pi^{2k}\times [-\infty,\infty]^{k(k-1)}$.
\eprop
{\bf Proof.} The proof is essentially the same as that for Proposition~\ref{P:Donsker2}. We proceed by induction.
Suppose that (\ref{conv3}) holds for a given $k\geq 2$. Let $z^\eps_{k+1}\in\Z^2_{\rm even}$ be such that $S_{\alpha, \sigma, \eps}z^\eps_{k+1} \to z_{k+1}$
for some $z_{k+1}=(x_{k+1}, t_{k+1})\in\R^2$. If $z_{k+1}=z_i$ for some $1\leq i\leq k$, then (\ref{conv3}) for $k+1$ follows from the induction assumption
and Proposition~\ref{P:Donsker2} applied to $C_{z^\eps_i}$ and $C_{z^\eps_{k+1}}$. Therefore we assume from now on $z_{k+1}\neq z_i$ for all $1\leq i\leq k$.

As in the proof of Proposition~\ref{P:Donsker2}, we construct $(C_{z^\eps_i})_{1\leq i\leq k}$ and $C_{z^\eps_{k+1}}$ from two independent percolation
configurations $\Omega^{\eps}_{[1]}$ and $\Omega^{\eps}_{[2]}$. First we construct $(C^{[1]}_{z^\eps_i})_{1\leq i\leq k}$ and $C^{[2]}_{z^\eps_{k+1}}$
respectively from $\Omega^{\eps}_{[1]}$ and $\Omega^{\eps}_{[2]}$, and then set $(C_{z^\eps_i})_{1\leq i\leq k}:=(C^{[1]}_{z^\eps_i})_{1\leq i\leq k}$. Next
we construct $C_{z^\eps_{k+1}}$ by successively exploring the status of edges in $\Omega^{\eps}_{[2]}$ until the first time we encounter an edge whose status in $\Omega^{\eps}_{[1]}$ has already been explored in the construction of $(C^{[1]}_{z^\eps_i})_{1\leq i\leq k}$, from which step onward, the exploration construction of $C_{z^\eps_{k+1}}$ will only use the status of edges that have already been explored, or if an edge is unexplored, then look up its status in $\Omega^{\eps}_{[1]}$. Let $\iota^{[12], \eps}$ be the first time $n$ when $C_{z^\eps_{k+1}}(n)$ intersects $(C_{z^\eps_i})_{1\leq i\leq k}$, then $C_{z^\eps_{k+1}}(\cdot) = C^{[2]}_{z^\eps_{k+1}}(\cdot)$ and $r_{z^\eps_{k+1}}(\cdot) = r^{[2]}_{z^\eps_{k+1}}(\cdot)$ on $[n^\eps_{k+1}, \iota^{[12],\eps}-1]$.

As in the proof of Proposition~\ref{P:Donsker2}, it suffices to go to a weakly convergent subsequence of $S_{\alpha, \sigma, \eps}\big((\gamma_{z^\eps_i}, r_{z^\eps_i})_{1\leq i\leq k+1}, (\kappa_{\gamma\gamma}^{\eps, ij}, \kappa_{rr}^{\eps, ij})_{1\leq i<j\leq k+1}\big)$ indexed by $(\eps_m)_{m\in\N}$ such that (\ref{errorBC}) and (\ref{gapBC}) hold with $\max_{i=1,2}$ therein replaced by $\max_{1\leq i\leq k+1}$. For such a sequence of $(\eps_m)_{m\in\N}$, we then apply Skorohod's
representation theorem to couple the sequence of percolation configurations $(\Omega^{\eps_m}_{[1]}, \Omega^{\eps_m}_{[2]})_{m\in\N}$, such that
\be\label{preskorohod2}
\begin{aligned}
& S_{\alpha, \sigma, \eps_m}\big((\gamma^{[1]}_{z^{\eps_m}_i}, r^{[1]}_{z^{\eps_m}_i})_{1\leq i\leq k}, (\kappa_{\gamma\gamma}^{[1],\eps_m, ij}, \kappa_{rr}^{[1],\eps_m, ij})_{1\leq i<j\leq k},\ (\gamma^{[2]}_{z^{\eps_m}_{k+1}}, r^{[2]}_{z^{\eps_m}_{k+1}})\big) \\
\longrightarrow &\quad  ((B^{[1]}_i, B^{[1]}_i)_{1\leq i\leq k}, (\kappa^{[1], ij})_{1\leq i<j\leq k}, (B^{[2]}_{k+1}, B^{[2]}_{k+1})) \qquad {\rm a.s.},
\end{aligned}
\ee
where $(B^{[1]}_i)_{1\leq i\leq k}$ is a collection of coalescing Brownian motions starting at $(z_i)_{1\leq i\leq k}$ with pairwise coalescence time $\kappa^{[1],ij}$,
and $B^{[2]}_{k+1}$ is an independent Brownian motion starting at $z_{k+1}$, all defined on the same probability space as $(\Omega^{\eps_m}_{[1]}, \Omega^{\eps_m}_{[2]})_{m\in\N}$. Let
$$
\kappa^{[12]} := \inf\{t\in\R : B^{[2]}_{k+1}(t)=B^{[1]}_i(t) \ \mbox{for some } 1\leq i\leq k\},
$$
and assume that $B^{[2]}_{k+1}(\kappa^{[12]}) = B^{[1]}_{i_0}(\kappa^{[12]})$ for some $1\leq i_0\leq k$. Then setting $(B_i)_{1\leq i\leq k}:= (B^{[1]}_i)_{1\leq i\leq k}$, $B_{k+1}(\cdot):= B^{[2]}_{k+1}(\cdot)$ on $[t_{k+1}, \kappa^{[12]}]$ and $B_{k+1}(\cdot):= B^{[1]}_{i_0}(\cdot)$ on $[\kappa^{[12]},\infty)$ produces a
collection of coalescing Brownian motions $(B_i)_{1\leq i\leq k+1}$ starting respectively at $(z_i)_{1\leq i\leq k+1}$. With $(C_{z^{\eps_m}_i})_{1\leq i\leq k+1}$
constructed from $(\Omega^{\eps_m}_{[1]}, \Omega^{\eps_m}_{[2]})$ as before and the coupling we have, all it remains is to prove that the convergence in (\ref{conv3}) takes place a.s.\ along the sequence indexed by $\eps_m$. By (\ref{preskorohod2}), it suffices to verify the a.s.\ convergence of
$$
S_{\alpha, \sigma, \eps_m}\big(\gamma_{z^{\eps_m}_{k+1}}, r_{z^{\eps_m}_{k+1}}, (\kappa^{\eps_m, i(k+1)}_{\gamma\gamma}, \kappa^{\eps_m, i(k+1)}_{rr})_{1\leq i\leq k}\big)
\longrightarrow (B_{k+1}, B_{k+1}, (\kappa^{i(k+1)}, \kappa^{i(k+1)})_{1\leq i\leq k}).
$$
By the a.s.\ convergence of coalescence times among $(\gamma_{z^{\eps_m}_i}, r_{z^{\eps_m}_i})_{1\leq i\leq k}=(\gamma^{[1]}_{z^{\eps_m}_i}, r^{[1]}_{z^{\eps_m}_i})_{1\leq i\leq k}$ in (\ref{preskorohod2}), the above a.s.\ convergence may be reduced further to showing the a.s.\ convergence of
$$
S_{\alpha, \sigma, \eps_m}\big(\gamma_{z^{\eps_m}_{k+1}}, r_{z^{\eps_m}_{k+1}}, \kappa^{\eps_m, i_0(k+1)}_{\gamma\gamma}, \kappa^{\eps_m, i_0(k+1)}_{rr})\longrightarrow
(B_{k+1}, B_{k+1}, \kappa^{i_0(k+1)}, \kappa^{i_0(k+1)}),
$$
which concerns only the pair of exploration clusters $C_{z^{\eps_m}_{i_0}}$ and $C_{z^{\eps_m}_{k+1}}$. As in the proof of Proposition~\ref{P:Donsker2}, a.s.\ either $y_{i_0}:= B_{i_0}(t_{i_0}\vee t_{k+1}) < y_{k+1}:=B_{k+1}(t_{i_0}\vee t_{k+1})$ or $y_{k+1}<y_{i_0}$. Treating the two cases separately,
the rest of the argument is then exactly the same as in the proof of Proposition~\ref{P:Donsker2}.
\qed

\section{Verification of (B1) and (B2)}\label{S:B12}
In this section, we conclude the proof of Theorem~\ref{T:main} by verifying conditions (B1) and (B2) in Theorem~\ref{T:convct1} for ${\cal X}_\eps:=S_{\alpha, \sigma,\eps}\overline{\Gamma}$ as $\eps\downarrow 0$.
\bigskip

\noindent
{\bf Verification of (B1).} In our setting, condition (B1) amounts to showing that for all $t>0$,
\be\label{B1}
\limsup_{\eps\downarrow 0} \sup_{(a,t_0)\in\R^2} \P(\eta_{{\cal X}_\eps}(t_0, t; a, a+\delta)\geq 2) \to 0 \quad \mbox{as}\quad \delta\downarrow 0.
\ee
Note that, if we denote $a_\eps:= \alpha t_0\eps^{-1} + a\sigma\eps^{-1/2}$ and $\delta_\eps:= \delta\sigma\eps^{-1/2}$, then
$$
\begin{aligned}
\eta_{{\cal X}_\eps}(t_0, t; a, a+\delta) & = \eta_{\Gamma}(t_0\eps^{-1}, t\eps^{-1}; a_\eps, a_\eps+\delta_\eps) \leq \eta_{\Gamma}(\lfloor t_0\eps^{-1}\rfloor, \lfloor t\eps^{-1}\rfloor; \lfloor a_\eps\rfloor, \lfloor a_\eps\rfloor +\lceil \delta_\eps\rceil+1),
\end{aligned}
$$
where we used the fact that paths in $\Gamma$ are nearest-neighbor paths and coalesce when they intersect. Therefore by translation invariance of $\Gamma$
under shifts by vectors in $\Z^2_{\rm even}$, uniformly in $(a, t_0)\in\R^2$, we have
\be\label{B1approx}
\begin{aligned}
\P(\eta_{{\cal X}_\eps}(t_0, t; a, a+\delta)\geq 2) & \leq \P(\eta_{\Gamma}(\lfloor t_0\eps^{-1}\rfloor, \lfloor t\eps^{-1}\rfloor; \lfloor a_\eps\rfloor, \lfloor a_\eps\rfloor +\lceil \delta_\eps\rceil+1)\geq 2) \\
& \leq \P(\eta_{\Gamma}(0, \lfloor t\eps^{-1}\rfloor; 0, \lceil \delta_\eps\rceil+2)\geq 2) \\
& \leq \P(\eta_{\cal R}(0, \lfloor t\eps^{-1}\rfloor; 0, \lceil \delta_\eps\rceil+2)\geq 2),
\end{aligned}
\ee
where ${\cal R}:=\{r_{(x,0)}: (x,0)\in\Z^2_{\rm even}\}$, and in the last inequality we applied Lemma~\ref{L:coalesce}. Assume that $x_\eps:=\lceil \delta_\eps\rceil+2$ is even, otherwise replace the constant $2$ by $3$. Note that $\eta_{\cal R}(0, \lfloor t\eps^{-1}\rfloor; 0, x_\eps)\geq 2$ if and only if $r_o$ and $r_{(x_\eps, 0)}$ do not coalesce before or at time $\lfloor t\eps^{-1}\rfloor$. Since $S_{\alpha, \sigma, \eps}(x_\eps, 0)\to (\delta,0)$ as $\eps\downarrow 0$, Prop.~\ref{P:Donsker2} implies that $\eps\kappa^\eps_{rr}$, where $\kappa^\eps_{rr}$ is the time of coalescence of $r_o$ and $r_{(x_\eps, 0)}$, converges in distribution to the time of coalescence $\kappa(\delta)$ of two coalescing Brownian motions starting respectively at $(0,0)$ and $(\delta, 0)$. Since $\eps \lfloor t\eps^{-1}\rfloor\to t>0$, we have
\be\label{B1kappa}
\lim_{\eps\downarrow 0}\P(\eta_{\cal R}(0, \lfloor t\eps^{-1}\rfloor; 0, x_\eps)\geq 2) = \P(\kappa(\delta)\geq t),
\ee
which tends to $0$ linearly in $\delta$ as $\delta\downarrow 0$; (\ref{B1}) then follows.
\qed
\bigskip

\noindent
{\bf Verification of (B2).} The key observation here is that, by replacing paths in $\Gamma$ by the associated exploration clusters, we end up with increasing and decreasing events of the underlying percolation configuration, for which we can then apply the FKG inequality. The details are as follows.

Condition (B2) in Theorem~\ref{T:convct1} amounts to showing that for all $t>0$,
\be\label{B2}
\delta^{-1} \limsup_{\eps\downarrow 0} \sup_{(a,t_0)\in\R^2} \P(\eta_{{\cal X}_\eps}(t_0, t; a, a+\delta)\geq 3) \to 0 \quad \mbox{as}\quad \delta\downarrow 0.
\ee
Let $a_\eps$, $\delta_\eps$, $x_\eps$ and $\cal R$ be as in the verification of (B1). Then similar to (\ref{B1approx}), we have
$$
\P(\eta_{{\cal X}_\eps}(t_0, t; a, a+\delta)\geq 3) \leq \P(\eta_{\cal R}(0, \lfloor t\eps^{-1}\rfloor; 0, x_\eps)\geq 3)
$$
uniformly in $(a, t_0)\in\R^2$. By (\ref{B1kappa}), to prove (\ref{B2}), it then suffices to show that for all $n\in\N$ and $x\in\N$,
\be\label{FKG}
\P(\eta_{\cal R}(0, n; 0, 2x)\geq 3) \leq \P(\eta_{\cal R}(0, n; 0, 2x)\geq 2)^2.
\ee
To simplify notation, let $C_i:=C_{(2i,0)}$, $l^n_i:=l^n_{(2i,0)}$, and $r_i:=r_{(2i, 0)}$ denote respectively the exploration cluster at $(2i,0)$ and its left and
right boundaries. Since paths in ${\cal R}$ are ordered by Lemma~\ref{L:coalesce}, we can write
$$
\P(\eta_{\cal R}(0, n; 0, 2x)\geq 3) = \sum_{k=1}^{x-1} \P(r_0(n)=r_{k-1}(n)<r_k(n)<r_x(n)).
$$
For each $1\leq k\leq x-1$, the event $\{r_0(n)=r_{k-1}(n)<r_k(n)<r_x(n)\}$ is the same as the event that the exploration clusters $C_0$ and $C_{k-1}$
intersect before or at time $n$, but $C_{k-1}$, $C_{k}$ and $C_{x}$ are mutually disjoint up to time $n$. Since different exploration
clusters evolve independently before they intersect, we can write
\begin{eqnarray}
&& \P(r_0(n)=r_{k-1}(n)<r_k(n)<r_x(n)) \nonumber \\
&=& \P(r_0(n)=r_{k-1}(n), \ r_{k-1}(\cdot)< l^{n}_k(\cdot) \mbox{ and } r_{k}(\cdot)< l^{n}_x(\cdot) \mbox{ on }
[0,n]) \nonumber \\
&=& \P(r^{[1]}_0(n)=r^{[1]}_{k-1}(n), \ r^{[1]}_{k-1}(\cdot)< l^{[2], n}_k(\cdot) \mbox{ and } r^{[2]}_{k}(\cdot)< l^{[3], n}_x(\cdot) \mbox{ on }
[0,n]),  \label{preFKG}
\end{eqnarray}
where $(C^{[1]}_0, C^{[1]}_{k-1})$, $C^{[2]}_k$ and $C^{[3]}_x$ and their boundaries are constructed on three independent percolation edge configurations $\Omega^{[1]}, \Omega^{[2]}$ and $\Omega^{[3]}$, defined as in (\ref{Omega}). Conditional on the realization of $\Omega^{[1]}$ and $\Omega^{[3]}$, and hence the realization of
$r^{[1]}_0, r^{[1]}_{k-1}$ and $l^{[3], n}_x$, we observe that the event $\{r^{[1]}_{k-1}(\cdot)< l^{[2], n}_k(\cdot) \mbox{ on } [0,n]\}$ is increasing in the edge configuration $\Omega^{[2]}$, while the event $\{r^{[2]}_{k}(\cdot)< l^{[3], n}_x(\cdot) \mbox{ on } [0,n]\}$ is decreasing in $\Omega^{[2]}$. Indeed, as more edges are switched from closed to open in $\Omega^{[2]}$, both $l^{[2],n}_k$
and $r^{[2]}_k$ can only increase. Therefore by the FKG inequality applied to the i.i.d.\ Bernoulli random variables in $\Omega^{[2]}$, and using the independence of $\Omega^{[1]}$ and $\Omega^{[3]}$, we have
$$
\begin{aligned}
&\ \P(r^{[1]}_0(n)=r^{[1]}_{k-1}(n), \ r^{[1]}_{k-1}(\cdot)< l^{[2], n}_k(\cdot) \mbox{ and } r^{[2]}_{k}(\cdot)< l^{[3], n}_x(\cdot) \mbox{ on }
[0,n]) \\
\leq &\ \E\big[1_{\{r^{[1]}_0(n)=r^{[1]}_{k-1}(n)\}}\ \P\big(r^{[1]}_{k-1}(\cdot)< l^{[2], n}_k(\cdot) \mbox{ on } [0,n]\,\big|\,\Omega^{[1]}\big) \ \P\big(r^{[2]}_{k}(\cdot)< l^{[3], n}_x(\cdot)\mbox{ on } [0,n] \,\big|\,\Omega^{[3]}\big) \big] \\
= &\ \P\big( r^{[1]}_0(n)=r^{[1]}_{k-1}(n), \mbox{ and } r^{[1]}_{k-1}(\cdot)< l^{[2], n}_k(\cdot) \mbox{ on } [0,n]\big)\ \P\big(r^{[2]}_{k}(\cdot)< l^{[3], n}_x(\cdot)\mbox{ on } [0,n]\big) \\
= &\ \P\big( r_0(n)=r_{k-1}(n), \mbox{ and } r_{k-1}(\cdot)< l^{n}_k(\cdot) \mbox{ on } [0,n]\big)\ \P\big(r_{k}(\cdot)< l^{n}_x(\cdot)\mbox{ on } [0,n]\big) \\
= &\ \P( r_0(n)=r_{k-1}(n)< r_k(n) )\, \P(r_k(n)<r_x(n)) \\
\leq &\ \P( r_0(n)=r_{k-1}(n)< r_k(n) )\, \P(r_0(n)<r_x(n)).
\end{aligned}
$$
Summing the above inequality over $1\leq k\leq x-1$ then gives (\ref{FKG}), and hence (\ref{B2}).
\qed

\bigskip

\noindent
{\bf Acknowledgement} A.~Sarkar thanks the Institute for Mathematical Sciences, National University of Singapore, for support during the workshop
{\em Probability and Discrete Mathematics in Mathematical Biology} in May 2011, when this project was initiated. R.~Sun is supported by grant R-146-000-119-133
from the National University of Singapore.

\end{document}